\documentclass[leqno]{aadmbook}
\usepackage{amsthm}
\usepackage{amsfonts}
\usepackage{amsmath}
\usepackage{latexsym,amsfonts} 
\usepackage{epsfig,float,graphicx}
\restylefloat{figure} 
\def\Dj{\rlap{--}D}
\def\ds{\displaystyle}
\def\dzn{,\kern-0.1em,}
\def\r{ \; \;  }

\newtheorem{thm}{Theorem}
\newtheorem{prp}{Proposition}
\newtheorem{clm}{Claim}
\newtheorem{lem}{Lemma}

\newtheorem{defn}{Definition}
\newtheorem{conj}{Conjecture}

\def\proof{\noindent\textbf{Proof}.\ \ }
\def\qed{\hfill$\Box$}

\input cyracc.def
 \font\eightcyr=wncyr10 scaled 833
\def\cyr{\eightcyr\cyracc} 

\def\cC{\mathcal{C}}
\def\cD{\mathcal{D}}
\def\cF{\mathcal{F}}
\def\cG{\mathcal{G}}
\def\cH{\mathcal{H}}

\def\cL{\mathcal{L}}
\def\cR{\mathcal{R}}
\def\cS{\mathcal{S}}
\def\cW{\mathcal{W}}

\def\be{\begin{equation} }
\def\ee{\end{equation} }
\def\bfl{\begin{flushleft} }
\def\efl{\end{flushleft} }
\def\bfr{\begin{flushright} }
\def\efr{\end{flushright} }
\def\bc{\begin{center}}
\def\vs*{\vspace*}
\def\hs*{\hspace*}
\def\ec{\end{center}}
\def\beq{\begin{eqnarray}}
\def\eeq{\end{eqnarray}}

\def\ben{\begin{enumerate}}
\def\een{\end{enumerate}}
\def\bit{\begin{itemize}}
\def\eit{\end{itemize}}

\begin{document}


\oddsidemargin 16.5mm
\evensidemargin 16.5mm

\thispagestyle{plain}

\begin{center}
{\large \sc  Applicable Analysis and Discrete Mathematics}
{\small available online at  http:/$\!$/pefmath.etf.rs }
\end{center}

\noindent{\small{\sc  Appl.\ Anal.\ Discrete Math.\ }{\bf x} (xxxx),
xxx--xxx.} \hfill{\scriptsize doi:10.2298/AADMxxxxxxxx}

\vspace{5cc}
\begin{center}
{\large\bf
  ENUMERATION OF HAMILTONIAN CYCLES \\
  ON A THICK GRID CYLINDER --- PART II: \\
  CONTRACTIBLE HAMILTONIAN CYCLES
\rule{0mm}{6mm}\renewcommand{\thefootnote}{}
\footnotetext{

\rule{2.4mm}{0mm}Keywords and Phrases: Hamiltonian cycles,
transfer matrix method, contractible cycle }}

\vspace{1cc}
{\large\it Olga Bodro\v{z}a-Panti\'{c}, Harris Kwong, \\ \vspace*{3pt}
 Jelena \Dj oki\' c,  Rade Doroslova\v{c}ki, Milan Panti\'{c}}

\vspace{1cc}
\parbox{24cc}{\small
In this series of papers, the primary goal is to enumerate Hamiltonian
cycles (HC's) on the grid cylinder graphs $P_{m+1}\times C_n$, where
$n$ is allowed to grow whilst $m$ is fixed.  In Part~I, we studied the
so-called non-contractible HC's.  Here, in Part~II, we proceed further
on to the contractible case.  We propose two different novel
characterizations of contractible HC's, from which we construct
digraphs for enumerating the contractible HC's.  Given the impression
which the computational data for $m \leq 9$ convey, we conjecture
that the asymptotic domination of the contractible HC's versus the
non-contractible HC's, among the total number of HC's, depends on the
parity of $m$.}
\end{center}

\vspace{1cc}



\vspace{1.0cc}
\begin{center}
{\bf 1. INTRODUCTION}
\end{center}

Determining and enumerating Hamiltonian cycles in some specific grid
graphs (such as thick grid cylinder graphs, which are studied here) is
of quite some relevance to statistical physics \cite{J1} and polymer
science \cite{BPPB}.  An ample amount of references related to this
topic may be found in Part~I~\cite{BKDP}.  A few novel applications of
this type of research can be found within the field of network
systems, which revolves around computer network functionality.
Hamiltonian cycles play a vital role there, because they cover all the
nodes of the system. In~\cite{PMBSJ} the issue of handling
indeterminacy for interval data under neutrosophic environment is
considered.  Another field, which may benefit from our research, is
that of cyber security.  There, digital microfluidic biochips (DMFBs)
are making the transition to the marketplace for commercial
exploitation.  For example, the microelectrode dot array (MEDA) is a
next-generation DMFB platform that supports real-time sensing of
droplets and has the added advantage of important security
protection~\cite{LCK}.

\begin{figure}[htb]
\begin{center}
\includegraphics[width=4.5in]{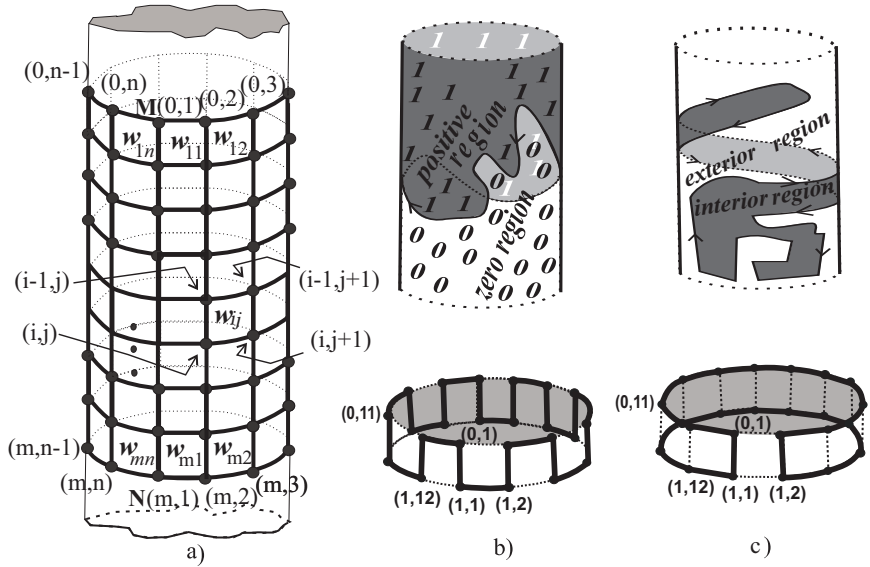}
\end{center}
\caption{(a) The graph $P_{m+1} \times C_n$ with its cells (windows)
  labeled by $w_{i,j}$, where $1\leq i\leq m$ and $1\leq j\leq n$.
  \quad (b) A non-contractible, closed Jordan curve on an infinite
  cylindrical surface.  \quad (c) A contractible, closed Jordan curve
  on an infinite cylindrical surface.}
\label{fig:2HC}
\end{figure}

When $m$ is fixed, the graphs $P_{m+1} \times C_n$ are referred to as
the \textbf{\emph{thick grid cylinders}} (see Figure~\ref{fig:2HC}a).
When $n\geq2$, there are two kinds of Hamiltonian cycles on such
graphs.  The first kind, denoted by HC$^{\; nc}$'s, are not
contractible when perceived as closed Jordan curves (see
Figure~\ref{fig:2HC}b) on the infinite cylindrical surface on which
the graph $P_{m+1} \times C_n$ is settled.  They were examined in
Part~I~\cite{BKDP} of this series.  The second kind of HC's, denoted
by HC$^{\; c}$'s, are the contractible ones.  They are studied in
Part~II of the series (this exposition).  In both parts, we study the
topological properties of the HC's.  Based on these properties, we
construct digraphs from which the HC's can be counted.  The motivation
behind our investigations is made clear in Part~I, together with the
reasons why we have opted for the cell-coding approach.

Contractible HC's are more complicated than the non-contractible ones.
These contractible HC's divide the underlying infinite cylindrical
surface into two separate regions.  The first is bounded and is called
the \textbf{\emph{interior}}, whereas the second one is called the
\textbf{\emph{exterior}} of the HC in question (see
Figure~\ref{fig:2HC}b-c). Moreover, we refer to these regions as the
\textbf{\emph{zero}} and \textbf{\emph{non-zero region}} depending on
whether a zero is assigned to the squares of the interior or the
exterior region.

The paper is organised as follows: in Section~2, we examine
HC$^{\;c}$'s whose non-zero region is the exterior.  Section~3 is
devoted to HC$^{\; c}$'s whose non-zero region is their interior.

The notations $h^{nc}_m(n)$ and $h^c_m(n)$ stand for the number of
HC$^{\; nc}$'s and HC$^{\;c}$'s, respectively.  Their respective
generating functions are
  \[ \cH^{nc}_m(x) \stackrel{\rm def}{=} \sum_{n\geq1}
  h^{nc}_m(n+1) x^n,
  \qquad\mbox{and}\qquad
  \cH^c_m(x) \stackrel{\rm def}{=} \sum_{n\geq1} h^c_m(n+1) x^n. \]
The overall number of HC's in the thick grid cylinder graph $P_{m+1}
\times C_n$ ($m \geq 1$, $n \geq 2$) is denoted by $h_m(n)$.  Clearly,
$h_m(n) = h^{nc}_m(n) + h^c_m(n)$ and its generating function
$\ds\cH_m(x) = \sum_{n\geq 1} h_m(n+1) x^n$ fulfills the equation
$\cH_m(x) = \cH^{nc}_m(x) + \cH^c_m(x)$.

The orientation of a HC$^{\;c}$ is determined in such a way that when
traversing alongside the considered HC its interior region is always
on the right-hand side (see Figure~\ref{fig:2HC}c). Recall an
assertion from Theorem 1 of \cite{BKP} concerning $h_m^c$:
  \be
  \label{thm:zeroh} h^c_m(n) = 0 \mbox{  \ if and only if $m$ is even
  and $n$ is odd.}
  \ee

Further, let us be reminded of a few additional definitions from Part
I~\cite{BKDP} needed hereinafter.  All the rest may as well be found
in the said paper, unless explicitly stated differently.

\begin{defn}
Given an integer word $d_1d_2\ldots d_m$, its \textbf{\emph{support}}
is defined as the ternary word $\bar{d}_1\bar{d}_2\ldots\bar{d}_m$,
where
  \[ \bar{d}_i = \begin{cases}
  \phantom{-}1 & if \r   d_i>0, \cr
  \phantom{-}0 & if \r  d_i=0, \cr
  -1 & if  \r  d_i<0. \cr \end{cases} \]
The \textbf{\emph{support}} of an integer matrix $[d_{i,j}]$ is
defined in a similar fashion.
\end{defn}

\begin{defn}
The factor $u$ of a word $v$ is called a \textbf{\emph{$b$-factor}} if
it is a block of consecutive letters all of which are equal to $b$.  A
$b$-factor of $v$ is said to be \textbf{\emph{maximal}} if it is not a
proper factor of another $b$-factor of $v$.
\end{defn}

Recall that the \textbf{\emph{window lattice graph}} $W_{m,n}$, whose
vertices are the square cells (or \textbf{\emph{windows}}) $w_{i,j}$
($1 \leq i \leq m, \; \; 1 \leq j \leq n$) of $P_{m+1} \times C_n$, is
isomorphic to $P_m \times C_n$.  For a HC$^{\;c}$, the interior
windows (marked with 0's as in Figure~\ref{fig:ext-flat}) form the
\textbf{\emph{interior tree (IT)}} in $W_{m,n}$.  Nonetheless, the
exterior windows form a forest of \textbf{\emph{exterior trees
(ET's)}}.  Note that only one ET from this forest contains exactly
one window in the first as well as in the last (the $m^{\mbox{th}}$)
row of $W_{m,n}$, called the \textbf{\emph{up}} and \textbf{\emph{down
root}}, respectively.  We call this particular ET the
\textbf{\emph{split tree (ST)}} of the HC in question.  Any other ET
different from the split tree contains either exactly one down root or
exactly one up root, but not both.  The ET's with a down root are
called the \textbf{\emph{down trees (DT's)}}, whereas the ET's with an
up root are referred to as the \textbf{\emph{up trees (UT's)}}.

\vspace{0.5cc}
\noindent\emph{Example 1.} \ %
For the purpose of illustration, take a look at the HC$^{\; c}$
depicted in Figures~\ref{fig:ext-flat} and~\ref{fig:int-flat} whose
split tree has the down root $w_{10,3}$, and the up root $w_{1,1}$. It
also has one ET with the down root $w_{10,9}$ (hence a DT), and one ET
with the up root $w_{1,7}$ (hence a UT); they are labeled by non-zero
integers in Figure~3.

\vspace{0.5cc}

Note that, it suffices to examine only those HC$^{\; c}$'s in
$P_{m+1}\times C_n$ whose split tree has $w_{11}$ for its up root. Let
the number of such HC$^{\; c}$'s be $\varphi^c_m(n-2)$, where $m\geq
1$, and $n\geq2$, and let the associated generating function be:
   \[ \Phi^c_m(x) \stackrel{\rm def}{=}
   \sum_{k\geq 0} \varphi^c_m(k)x^k. \]
This implies that the total number of HC$^{\; c}$'s in $P_{m+1}\times
C_n$ is given by
$$h^c_m(n) = n\varphi^c_m(n-2).$$  Consequently,
$$ \ds  \cH_m^c(x) = \sum_{n\geq 1} h^c_m(n+1)x^n
  =  \sum_{n\geq 1} (n+1)\varphi^c_m(n-1) x^n  $$
  \[
  =   \frac{d}{dx} \sum_{n\geq 1} \varphi^c_m(n-1) x^{n+1}
  = \frac{d}{dx} \left(x^2\,\Phi^c_m(x)\right). \]

There are two possible ways in which we code (or label) the windows
with appropriate integers. The first, described in Section 2, is the
one in which the windows of the IT are labeled with zeros, whilst the
remaining windows are labeled with non-zero numbers. The second, which
we deal with in Section 3, is the one in which the zero windows belong
to the ET's, whereas the non-zero windows belong to the IT.  This way,
any HC$^{c}$ can be viewed as a sequence of $n$ columns comprising the
coded windows.  This sets up a one-to-one correspondence between the
set of HC$^{\; c}$'s and the set of sequences of $n$ labeled columns.

Recall from~\cite{BKDP} that ${\cal G}_{m}$ represents an infinite
grid graph with vertices from the set $\{ (i,j) \in \mathbb{Z}^2 \mid
\; 0 \leq j \leq m \}$, in which the square cell determined by the
points: $(j-1+ kn, m-i)$, $(j+ kn, m-i)$, $(j+ kn, m-i+1)$ and $(j-1+
kn, m-i+1)$, with $1 \leq i \leq m \mbox{ and } 1 \leq j \leq n$, is
labeled $w_{ij}^k$ and is called a \textbf{\emph{window}}, too. We
also say that $w_{ij}^k$ belongs to the $(j+nk)^{\mbox{th}}$ column of
$\cG_m$. The set $\{w_{ij}^k \mid i,j,k \in \mathbb{Z} \; \wedge \; 1
\leq i \leq m \; \wedge \; 1 \leq j \leq n \}$ presents the set of
vertices of another infinite grid graph denoted by ${\cal W}_{m}$.

\begin{figure}[htbp]
\begin{center}
\includegraphics[width=4.8in]{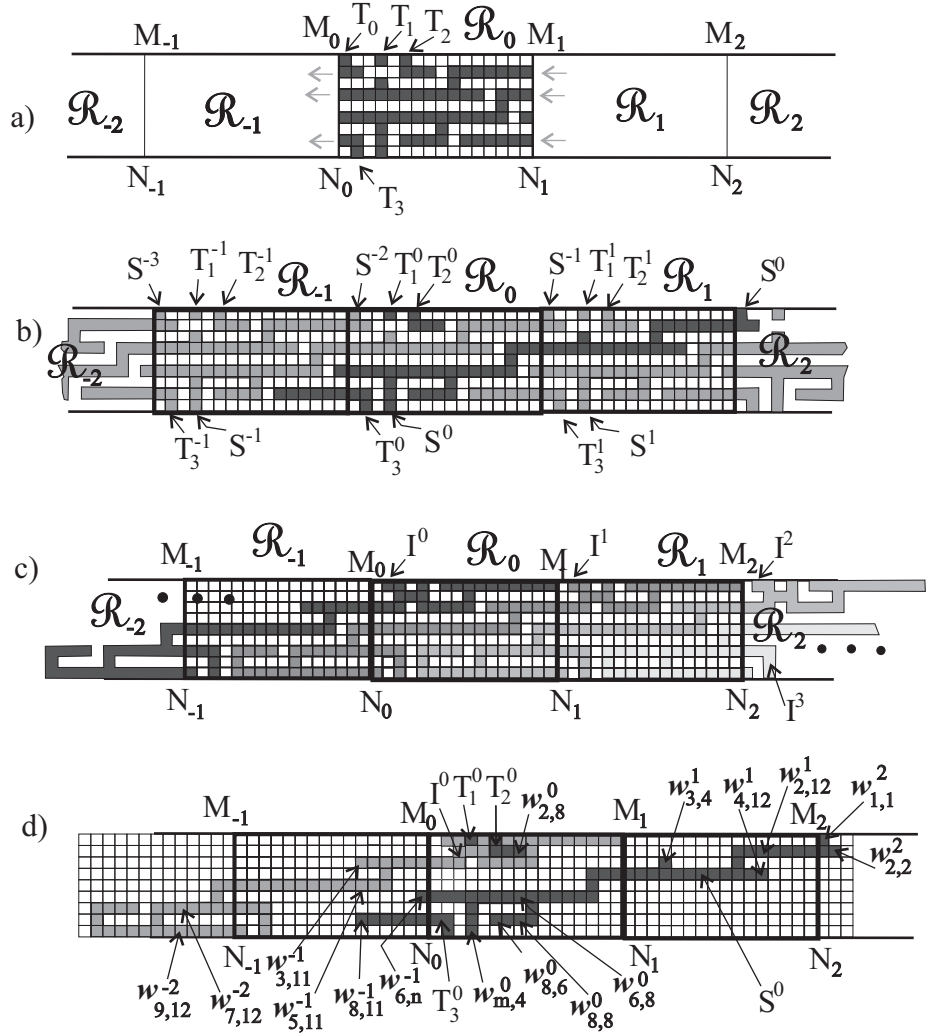}
\end{center}
\caption{%
(a) Unrolling and flattening the cylindrical surface that contains
a HC$^{\; c}$. \quad
(b) The rolling imprint of a HC$^{ \; c}$, with the copies of the
``new'' split tree and ETs in gray. \quad
(c) The rolling imprint of  a HC$^{\; c}$, with the copies of the
``new'' interior tree in gray. \quad
(d) The basis of a rolling imprint (BRI) consists of the windows of all
the ``new'' trees $I^0, S^0, T_1^0, T_2^0$ and $T_3^0$.}
\label{fig:RI}
\end{figure}

\begin{figure}[htbp]
\begin{center}
\includegraphics[width=4.8in]{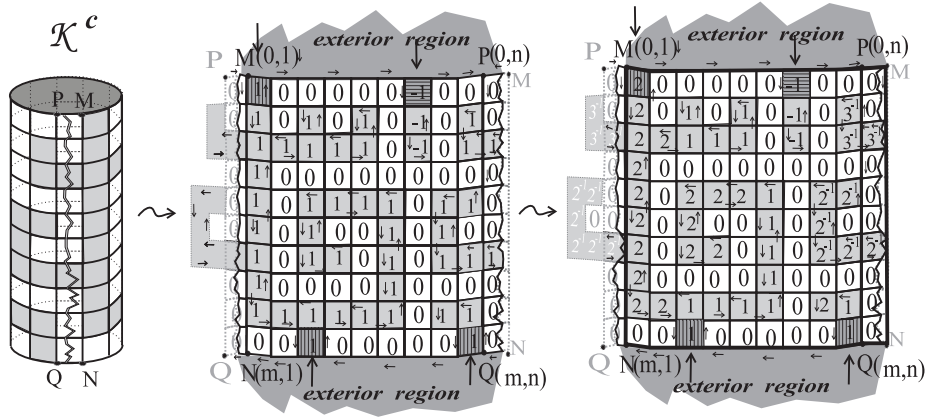}
\end{center}
\caption{A contractible Hamiltonian cycle of $P_{11}\times C_{10}$
  with the entries of the matrices  $A^{ \; c,Ext}$ (the first phase)
  and $B^{ \; c,Ext}$ (the second phase) written on their windows as
  the cylindrical surface is drawn on a flat surface.}
\label{fig:ext-flat}
\end{figure}

\begin{figure}[htbp]
\begin{center}
\includegraphics[width=4.8in]{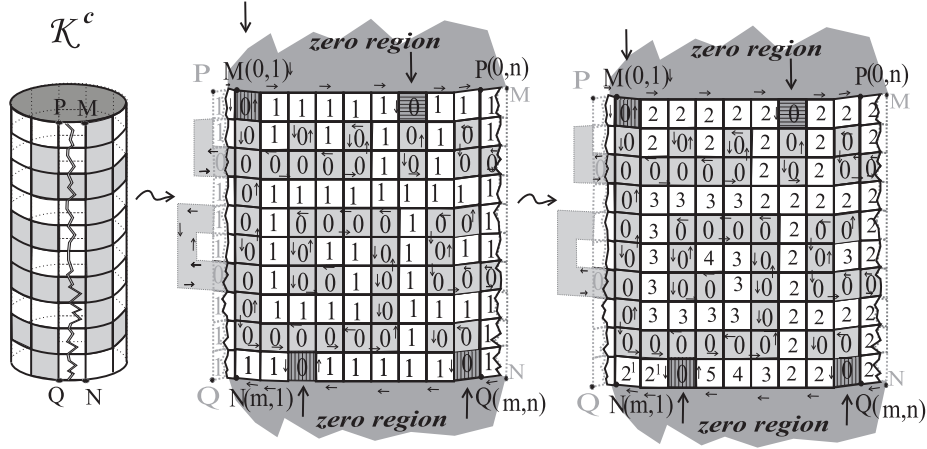}
\end{center}
\caption{A contractible HC of $P_{11}\times C_{10}$ with the entries
  of $A^{ \; c,Int}$ (the first phase) and $B^{ \; c,Int}$ (the second
  phase) inscribed on its windows, once the cylindrical surface was
  represented on a flat surface.}
\label{fig:int-flat},
\end{figure}

Consider a HC$^{\; c}$ in the graph $P_{m+1}\times C_n$.  Loosely
speaking, a \textbf{\emph{rolling imprint (RI)}} is a picture obtained
as follows.  First we ``cut through'' the surface of our graph
$P_{m+1} \times C_n$ (with a HC in it) along the line which connects
the vertices $M(0,1)$ and $N(m,1)$, see Figure~\ref{fig:2HC}a.  Next,
we unroll and flatten it; see the rectangle ${\cal R}_0: M_0N_0N_1M_1$
in Figure~\ref{fig:RI}a.  Finally we produce many copies of the
initial picture (${\cal R}_{-1}$, ${\cal R}_{-2}$, ${\cal R}_{-3},
\ldots$ and ${\cal R}_{1}$, ${\cal R}_{2}$, ${\cal R}_{3}, \ldots $),
and line them up to the left and to the right side accordingly; see
Figure~\ref{fig:RI}b.  Since the HC is contractible, its RI is
actually the graph ${\cal G}_{m}$ with infinitely many mutually
congruent polygonal lines on it.  These polygonal lines are the
boundaries of the polygons consisting of all the vertices
($w_{ij}^k$\,) of ${\cal W}_{m}$ that correspond to the windows
($w_{ij}$) of $P_{m+1}\times C_n$ from the interior of the HC$^{\; c}$
(the white squares in Figure~\ref{fig:RI}b or the gray squares in
Figure~\ref{fig:RI}c).  That way parts of the interior and exterior
trees that were initially broken by the process of ``cutting'' of the
$W_{m,n}$ are now assembled again into the original forms, and
multiplied in ${\cal W}_{m}$.  What we obtain is a sequence of copies
$ \ldots S^{-3}, S^{-2}, S^{-1}, S^{0}, S^{1}, S^{2}, \ldots $ of a
``new'' split tree $ S^{0}$, some sequences of copies $ \ldots
T_{s}^{-3}, T_{s}^{-2},T_{s}^{-1},T_{s}^{0},T_{s}^{1}, T_{s}^{2},
\ldots $ of the ``new'' exterior trees $ \ldots T_{s}^{0} $ and a
sequence of copies $ \ldots I^{-2},I^{-1}, I^{0}, I^{1}, I^{2}, \ldots
$ of the ``new'' interior tree $I^{0}$.

At this point we need to modify a few definitions stated
in~\cite{BKDP}, as follows:

\begin{defn}
The \textbf{\emph{ basis of a rolling imprint (BRI) }} is the union of
the vertex set of the split tree whose down root is in $\cR_0$, the
vertex sets of all the exterior trees (different from the split tree)
each of which has its root in $\cR_0$, and the vertex set of the
interior tree whose leftmost window from the first row belongs to
${\cal R}_0$.
\end{defn}

Note that, in this way, we establish a bijection between the set of
vertices in $V(W_{m,n})$ and the BRI (see Figure~\ref{fig:RI}d).

The aforementioned coding of the windows is, in both cases, dealt with
in two stages. In the first stage the graph $W_{m,n}$ is associated
with the matrix $A^{ c} = [a_{ij}]_{m \times n}$ whose entries are
from $\{-1,0,1\}$.  The windows $w_{i,j}$ are called the
\textbf{\emph{zero windows}} if and only if $a_{i,j} =0$, otherwise
they are named the \textbf{\emph{non-zero windows}}.  The coding is
done by associating the same number to each of the vertices of the
same tree (be it a ST, ET or IT) during the first stage. For instance,
in the first case all the vertices of the UT's were coded with $-1$,
whilst all the vertices of the DT's and ST were coded with
$1$. Therefore, we say that the DT's and the ST are
\textbf{\emph{positive trees (PT)}}, whereas for the UT's we say that
they are \textbf{\emph{negative trees (NT)}}, or simply
\textbf{\emph{non-zero trees}}, irrespective of the case.  In the
second case, the term \textbf{\emph{positive tree (PT)}} or
\textbf{\emph{non-zero tree}} refers to the IT.  The term
\textbf{\emph{zero tree}} is used in a similar manner.  The roll
number depends on the type of cell (zero or non-zero).

\begin{defn}
The \textbf{\emph{roll number}} (or simply \textbf{\emph{roll}}) of a
window $w_{i,j} \in V(W_{m,n})$, denoted by $r(w_{i,j})$ (or simply
$r$ if the window is clear from the context) is a unique integer $k$
for which $w_{ij}^k$ belongs to a non-zero tree of the BRI, or in case
$w_{ij}$ is a zero window, we set $r=0$.  We shall also say that the
window $w_{i,j}$ \textbf{\emph{belongs to roll~$r$}}.
\end{defn}

\vspace{0.5cc}
\noindent\emph{Example 2.} \ %
For the HC$^{\; c}$ whose BRI is presented in Figure~\ref{fig:RI}d,
the roll numbers of some specific windows are summarized below.
  \[ \renewcommand{\arraycolsep}{2pt}
  \begin{array}{|*{12}{c|}} \hline
  \mbox{non-zero} & \mbox{coding}
  & \multicolumn{10}{c|}{\mbox{roll number}} \\ \cline{3-12}
  \mbox{tree(s)}  & \mbox{method}
  & w_{9,12} & w_{7,12} & w_{3,11} & w_{5,11} & w_{1,2} & w_{1,n}
  & w_{8,11} & w_{8,8}  & w_{4,12} & w_{2,2} \\ \hline
  \mbox{IT}  & \mbox{second}
  & -2 & -2 & -1 & -1
  & 0 & 0 & \phantom{-}0 & 0 & 0 & 0 \\ \hline
  \mbox{ETs} & \mbox{first}
  & \phantom{-}0 & \phantom{-}0 & \phantom{-}0 & \phantom{-}0
  & 0 & 0 & -1 & 0 & 1 & 2 \\ \hline
  \end{array} \]

\vspace{0.25cc}

\begin{defn}
Two non-zero vertices $w_{i,t}$ and $w_{j,s}$ of $W_{m,n}$ with
$a_{i,t} = a_{j,s}$ are said to be \textbf{\emph{joined at the $k$-th
column with the roll number $r$}}, or simply
\textbf{\emph{$k^r$-joined}}, where $1 \leq k \leq n$, and
$-\big\lfloor\frac{m}{2}\big\rfloor \leq r \leq \big\lfloor\frac{m}{2}
\big\rfloor$, if and only if their corresponding windows in the BRI
belong to the same component in the subgraph of $\cW_m$ induced by the
set of all non-zero windows $w^z_{x,y}$ from the BRI that satisfy both
$a_{x,y} = a_{i,t} = a_{j,s}$, and either (i) $z=r(w_{xy})<r$, or (ii)
$z=r(w_{xy})=r$ and $y\leq k$.
\end{defn}

\vspace{0.5cc}
\noindent\emph{Example 3.} \ %
Let us once again take a look at Figure~\ref{fig:RI}d assuming the
first way of coding ($a_{1,1} \neq 1$).  There, windows $w_{6,8}$ and
$w_{8,8}$ are not $8^{0}$-joined, but instead are $9^{0}$-joined.
Also, $w_{2,12}$ is  $12^{1}$-joined with $w_{4,12}$.

\vspace{0.5cc}

In Sections~2 and~3, as we said, we present two different
characterizations of HC$^{\; c}$ where $w_{11}$ is the up root of the
split tree. Both of them allow for the use of the transfer matrix
method with a view to obtaining the values of $h^{ \; c}_m(n)$'s.  In
Section~4, we determine the upper bound of the so-called color$^r$
words which appear in these procedures.  Sections~5 and~6 contain
comparative analysis of the numerical results obtained by using these
two characterizations and some other conclusions including two new
conjectures. Section~7 is devoted to closing remarks.


\vspace{1.0cc}
\begin{center}
{\bf 2. CODING THE EXTERIOR TREES BY NON-ZERO ENTRIES}
\end{center}

\vspace{0.5cc}
\begin{center}
{\bf 2.1.  The First Phase --- the Matrix $A^{ \; c,Ext}$}
\end{center}

For any integer $m\geq1$, we associate with each HC$^{\; c}$ in
$P_{m+1}\times C_{n}$ with $w_{1,1}$ as the up root of the split tree
a matrix $A^{ \; c,Ext} = [a_{ij}]_{m \times n}$ whose entries are
defined in the following way:
  \[ a_{i,j} \stackrel{\rm def}{=} \begin{cases}
  \hfill 0 & \mbox{ if $w_{i,j}$ belongs to the IT}, \cr
        -1 & \mbox{ if $w_{i,j}$ belongs to a UT}, \cr
  \hfill 1 & \mbox{ if $w_{i,j}$ belongs to a DT or the ST.} \cr
  \end{cases} \]
Obviously, $a_{1,1} \stackrel{\rm def}{=}1$, $a_{1,2} \stackrel{\rm
def}{=}0, a_{1,n} \stackrel{\rm def}{=}0$, and $a_{2,1} \stackrel{\rm
def}{=}1$.  Note that $w_{1,1}$ is the only positive window in the
first row (on the ``negative coast''). We adopt the convention that
$a_{i,n+1} \stackrel{\rm def}{=} a_{i,1}$, and $a_{i,0} \stackrel{\rm
def}{=} a_{i,n}$, for $1\leq i\leq m$.

\begin{lem}  \label{lem:firstC}
Every $ HC^{ \; c}$ on the thick grid cylinder graph $P_{m+1} \times
C_n$ (with $w_{1,1}$ as the up root of the split tree) determines a
matrix $A^{ \; c,Ext} = [a_{ij}]_{m \times n}$, with entries from the
set $\{ -1, 0, 1\}$, which satisfies the conditions below.
\begin{enumerate}
\item \textbf{First and Last Row Conditions} (FL$^{ \; c,Ext}$):
  \begin{enumerate}
  \item $a_{1,1}=1$, and $a_{1,2}=a_{1,n}=0$.
  \item For $2 \leq j \leq n$,
    \[ a_{1,j}\in\{0,-1\}, \quad\mbox{and}\quad
    (a_{1,j},a_{1,j+1}) \not\equiv (-1,-1). \]
  \item For  $1 \leq j \leq n$,
    \[ a_{m,j}\in\{0,1\}, \quad\mbox{and}\quad
    (a_{m,j},a_{m,j+1}) \not\equiv (1,1). \]
  \end{enumerate}
\item \textbf{Adjacency of Column Conditions} (AC$^{\;  c,Ext}$):
  \begin{enumerate}
  \item For $1 \leq i \leq m-1$ and $1 \leq j \leq n$,
    \begin{eqnarray*}
    \lefteqn{(|a_{i,j}|,|a_{i+1,j}|,|a_{i,j+1}|,|a_{i+1,j+1}|)} \\
    &\notin& \{(1,1,1,1),(0,0,0,0),(0,1,1,0),(1,0,0,1)\}.
    \end{eqnarray*}
  \item For $1 \leq i \leq m$ and $1 \leq j \leq n$,
    \[ a_{i,j} \cdot a_{i,j+1} \neq -1. \]
  \item For $1 \leq i \leq m-1$ and $1 \leq j \leq n$,
    \[ a_{i,j} \cdot  a_{i+1,j} \neq -1. \]
  \end{enumerate}
\item \textbf{Root Conditions} (RC$^{ \; c,Ext}$):
  \begin{enumerate}
  \item Each connected component of the subgraph of the graph
    $W_{m,n}$ induced by the windows corresponding to the non-zero
    entries ($\pm1$) of the matrix $A^{ \; c,Ext}$ is a tree.
  \item There exists exactly one such tree (the split tree) containing
    exactly one window from the first row ($a_{11}=1$), and exactly
    one window from the last row of $W_{m,n}$.
  \item Each of the remaining trees, if any such exist, contains
    exactly one window from either the first or the last row of
    $W_{m,n}$.
  \end{enumerate}
\end{enumerate}
Conversely, every matrix $[a_{ij}]_{m \times n}$  with entries from
the set $\{ -1, 0, 1\}$  which fulfills the  conditions FL$^{ \;
  c,Ext}$, AC$^{ \; c,Ext}$ and RC$^{ \; c,Ext}$,  determines a unique
HC$^{\; c}$ on the graph $P_{m+1}\times C_{n}$ whose split tree
contains the window $w_{1,1}$.
\end{lem}

\proof
The necessity of all three imposed conditions is easily verifiable and
is thus left to the reader.  Therefore, we move on to the proof of
their sufficiency.  Let us observe all the regions determined by all
the non-zero windows including the two half-cylinders from the one
side (at this moment, we cannot assume that there exists a unique such
region).  The first two conditions are local conditions ensuring that
the boundary of the said regions (the edges which belong to both
non-zero window or the boundary of one of the two half-cylinders, and
a zero window) determines a unique spanning $2$-regular subgraph of
$P_{m+1} \times C_n$, that is, a union of cycles.

The proof that this graph consists of only one component (consequently
establishing the uniqueness of both zero and non-zero regions) can be
derived constructively.  The case $n=2$ is trivial, so we can assume
that $n\geq 3$.  The condition AC$^{ \; c,Ext}$ implies that each of
the components of the subgraph of $W_{m,n}$ induced by the windows
corresponding to the non-zero entries consists only of either
$1$-windows or $(-1)$-windows.  Thus, we have justified the existence
of both positive (DT's and ST) and negative trees (UT's).

Let $w_{m,p}$ be the down root of the unique positive tree $T_0$ (the
split tree) with a window $w_{11}$ in the first row in it and $T_1,
T_2, \ldots , T_k$ be all the NT's (if any such tree exists at all)
with the up roots $w_{1,j_1}, w_{1,j_2}, \ldots ,w_{1,j_k}$,
respectively, for which $3 \leq j_1 < j_2 < \ldots < j_k \leq n$.  Let
$T'_1, T'_2, \ldots, T'_l$ be all the PT's different from $T_0$ (if
any such tree exists in the first place) with the down roots
$w_{m,i_1}, w_{m,i_2}, \ldots, w_{m,i_l}$, respectively, for which
$i_s \equiv i'_s$ (mod $n$), where $1 \leq s \leq l$ and $p+2 \leq
i'_1 < i'_2 < \ldots < i'_l \leq n +p -2$. Our task is to obtain the
unique curve (the broken line) which separates the regions of the two
kinds of windows (the zero and non-zero ones).

We can start from the point $M:(0,1)$ (the upper-left point of
the up root $w_{1,1}$ of $T_0$) and move to the lower-left point of
$w_{m,p}$ using the edges of $P_{m+1} \times C_n$ that belong to the
boundary of $T_0$.  From there we continue towards the point
$(m,i_l+1)$, and then visit all the vertices on the boundary of $T'_l$
finishing at the point $(m,i_l)$.  Next, we visit the boundary of
$T'_{l-1}$, $T'_{l -2}$, $\ldots$.  After having visited the tree
$T'_1$ we end up at the point $(m,i_1)$.  Then, we move further to the
point $(m,p+1)$ and continue towards the point $(0,2)$ using the
remaining edges of the boundary of ST.  From there we similarly
continue visiting the boundaries of trees $T_1, T_2, \ldots , T_k$,
respectively, ending up at the point $M$ again (see
Figure~\ref{fig:ext-flat}).  By doing so, we pass through all the
edges on the boundary of these regions, obtaining a contractible HC.
\qed

\vspace{0.5cc}
Recall that, for fixed values of $k$ and $r$ ($1 \leq k \leq n$, and
$-\big\lfloor\frac{m}{2}\big\rfloor \leq r \leq \big\lfloor\frac{m}{2}
\big\rfloor$), the relation $k^r$-joined represents an equivalence
relation on the set of all non-zero windows $w_{x,y}$ that satisfy
either (i) $r(w_{x,y}) < r$, or (ii) $r(w_{x,y}) = r$ and $y\leq k$
(that is, whose window from the BRI belongs to the
$(y+nk)^{\mbox{th}}$ column of $\cG_m$ or to the left of it).
Furthermore, every equivalence class belongs to exactly one ET.  Hence
if this equivalence class belongs to a PT its windows can be
$k^r$-joined with at most one down root.  If it belongs to an NT its
windows can be $k^r$-joined with at most one up (negative) root.
Further, because an ST is a PT, we treat its down root as its main
root and this is what we shall assume below.  Note that the roll
number of $w_{1,1}$ could be different from~0 (for example, the roll
of $w_{1,1}$ in Figure~\ref{fig:RI}d is $2$). But, all the other roots
of the ET's have their roll number equal to $0$.


\vspace{0.5cc}
\begin{center}
{\bf 2.2. The First Characterization of HC$^{\; c}$ with $w_{11}$ \\
as the Up Root of the Split Tree}
\end{center}

Let $\cC^+ \stackrel{\rm def}{=} \big\{2,3,\ldots, \big\lfloor
\frac{m}{2} \big\rfloor+1\}$ and $\cC^- \stackrel{\rm def}{=}
\big\{-2,-3,\ldots,-\big\lfloor \frac{m}{2} \big\rfloor-1\}$.  For
each HC$^{\; c}$ with the window $w_{1,1}$ belonging to the split
tree, we associate the matrix $A^{ \; c,Ext} = [a_{i,j}]_{m \times n}$
with the matrix $B^{ \; c,Ext} = [(b_{i,j}, r_{i,j})]_{m \times n}$,
where $b_{i,j} \in \cC^+\cup\cC^-\cup \{1,0,-1\}$ and $- \big\lfloor
\frac{m}{2} \big\rfloor \leq r_{i,j} \leq \big\lfloor \frac{m}{2}
\big\rfloor$. The former of the two satisfies the conditions FL$^{ \;
  c,Ext}$, AC$^{ \; c,Ext}$, and RC$^{ \; c,Ext}$, whereas the latter
is constructed in the following way:
\begin{enumerate}
\item Define $r_{i,j} = r(w_{i,j})$.
\item Set  $b_{i,j} = a_{i,j}=0$ if $w_{i,j}$ belongs to the IT.
\item If $w_{i,j}$, where   $w_{i,j} \neq w_{1,1}$, is the (up) root
  of an NT (that is, $i=1$ and $a_{i,j} = -1$) or the down root of a
  PT (that is, $i=m$ and $a_{i,j} = 1$),  set $b_{i,j} = a_{i,j}$. If
  $w_{i,j}$ is neither the down root of a PT  nor  the up root of an
  NT, but it is $j^r$-joined with such a root, where $r = r(w_{i,j})$,
  set  $b_{i,j} =a_{i,j}$.
\item For each fixed column, say column $j$:
  \begin{enumerate}
  \item Scan the remaining positive windows $w_{i,j}$ with the same
    roll number from bottom to top (that is, from $i=m$ to $i=1$), and
    set $b_{i,j}$ to $z+1$, where $z$ is the ordinal number of the
    $j^r$-joined equivalence class, $r = r(w_{i,j})$, to which
    it belongs to (hence, the labels of the $b_{i,j}$'s start from $2$).
  \item Scan the remaining negative windows $w_{i,j}$ with the same
    roll number, from top to bottom (from $i=1$ to $i=m$), and set
    $b_{i,j}$ to $z-1$, where $z$ is the negative value of the ordinal
    number of the $j^r$-joined equivalence class,
    $r=r(w_{i,j})$, to which it belongs to (hence, the labels of the
    $b_{i,j}$'s start from $-2$).
  \end{enumerate}
\end{enumerate}

\vspace{0.5cc}
\noindent\emph{Example 4.} \ %
In Figure~\ref{fig:ext-flat}, the entries in the matrix $B^{ \;
c,Ext}$ are written on their respective windows ($\ds
b_{ij}^{r_{ij}}$ stands in place of $(b_{ij},r_{ij})$, or just $\ds
b_{ij}$ if $r_{ij}=0$). Note that in the $9^{\mbox{th}}$ column there
exist three parts of the same PT (it is the split tree) with the same
roll number $-1$, but the windows of only two of them are
$9^{-1}$-joined.  Consequently, the same value is associated to their
entries in the matrix $B^{ \; c,Ext}$ ($b_{5,9}=b_{7,9}=2$, and
$b_{3,9}=b_{2,9}=3$).  Another example is shown in
Figure~\ref{fig:RI}d. There, the entries of the matrix $B^{ \; c,Ext}
= [(b_{i,j}, r_{i,j})]_{9 \times 16}$ for windows $w_{1,1}$ and
$w_{2,2}$ are $(1,2)$, for windows $w_{3,4}$, $w_{2,12}$ and
$w_{4,12}$ are $(1,1)$, for windows $w_{m,2}$, $w_{m,4}$ and $w_{6,8}$
are $(1,0)$, for windows $w_{8,6}$ and $w_{8,8}$ are $(2,1)$, for
`windows $w_{1,4}$, $w_{1,6}$ and $w_{2,8}$ are $(-1,0)$ and for
windows $w_{3,11}$, $w_{5,11}$, $w_{7,12}$ and $w_{9,12}$ are $(0,0)$.

\vspace{0.5cc}

Consider all the existing maximal $b$-factors, where $b>0$ ($b<0$), in
the $j^{\mbox{th}}$ column $v = b_{1,j}b_{2,j}\ldots b_{m,j}$, where
$1\leq j\leq n$, of the matrix $[b_{i,j}]_{m \times n}$ corresponding
to the matrix $B^{ \; c,Ext} = [(b_{i,j},r_{i,j})]_{m \times n}$.  Let
them be, in their order of appearance (that is, from bottom to top for
positive windows, but from top to bottom for negative windows),
$p_1$-factor, $p_2$-factor, \ldots, $p_k$-factor, where $k\geq1$, and
$p_i\geq1$ ($p_i\leq -1$) for each $i$.  In addition, let
$r_1,r_2,\ldots,r_k$ denote the roll numbers associated with these
maximal factors.  The words $p_1 p_2 \ldots p_k$ and $r_1 r_2 \ldots
r_k$ are called the \textbf{\emph{positive (respectively, negative)
truncated word}} and the \textbf{\emph{positive (resp., negative)
truncated roll word}}, respectively, corresponding to the
$j^{\mbox{th}}$ column of $B^{ \; c,Ext}$. A subsequence of a
truncated word induced by the letters with the same roll number $r$ is
called a \textbf{\emph{positive (resp., negative) color$^r$ word}}.

\vspace{0.5cc}
\noindent\emph{Example 5.} \ %
For the $9^{\mbox{th}}$ column in Figure~\ref{fig:ext-flat}, the
positive truncated word, the positive truncated roll word, the
positive color$^0$ word, and the positive color$^{-1}$ word are 1223,
$0\,-1\,-1\,-1$, $1$ and 223, respectively.  Note that, in general,
$r_{1,1}$ need not be 0, and $b_{1,1}$ need not be 2.

\vspace{0.5cc}
\begin{center}
{\bf 2.3. Properties of the Matrix $B^{ \; c,Ext}$}
\end{center}

From the definition of the matrix $B^{ \; c,Ext} =
[(b_{i,j},r_{i,j})]_{m \times n}$, we can easily obtain a number of
properties expressed in the following theorem. Bear in mind that here
$(b_{i,n+1},r_{i,n+1}) \stackrel{\rm def}{=} (b_{i,1},r_{i,1})$, and
$(b_{i,0},r_{i,0}) \stackrel{\rm def}{=} (b_{i,n}, r_{i,n})$.

\begin{thm}  \label{thm:secondC}
The matrix $B^{ \; c,Ext} = [(b_{i,j},r_{i,j})]_{m \times n}$
satisfies the following conditions.
\begin{enumerate}

\item \textbf{Basic Properties}  
\begin{enumerate}
\item  
  The support of the matrix $[b_{i,j}]_{m \times n}$, that is, the matrix
  $[a_{i,j}]_{m \times n}$, satisfies the conditions FL$^{ \; c,Ext}$ and
  AC$^{ \; c,Ext}$.
  \item  
  \emph{Harmonization of the adjacent entries which have the same sign}:
  For $2\leq i\leq  m$, and $1\leq k\leq n$, if $a_{i-1,k} = a_{i,k}$,
  then $(b_{i-1,k},r_{i-1,k}) = (b_{i,k},r_{i,k})$.
 \item  
  For $1\leq i\leq m$, and $1\leq j\leq n$, if $a_{i,j} = 0$, then
  $r_{i,j} = 0$.
  \item  
   For  $3\leq j\leq (n-1)$, if $a_{1,j} \neq 0$, then
  $(b_{1,j},r_{1,j}) = (-1,0)$. \\
  For  $1\leq j\leq n$, if $a_{m,j} \neq 0$, then
  $(b_{m,j},r_{m,j})= (1,0)$.
\end{enumerate}  

\item \textbf{Column Properties} \par 
For $1\leq k\leq n$, the $k$-th column $[(b_{1,k},r_{1,k}),
(b_{2,k}, r_{2,k}), \ldots, (b_{m,k},r_{m,k})]^T$ of the matrix
$B^{ \; c,Ext}$ satisfies these conditions:
\begin{enumerate}
\item   
  If there exists an entry $(s,r)$ in the $k^{\mbox{th}}$ column of
  the matrix $B^{ \; c,Ext}$, where $s\geq3$, then for each $\ell\in
  \{2,3,\ldots,s-1\}$, at least one copy of the entry $(\ell,r)$ must
  appear after the last appearance of the entry $(s,r)$.  Likewise, if
  there exists an entry $(s,r)$ in the $k^{\mbox{th}}$ column of the
  matrix $B^{ \; c,Ext}$, where $s\leq-3$, then for each
  $\ell\in\{-2,-3,\ldots,s+1\}$, at least one copy of the entry
  $(\ell,r)$ must appear before the first appearance of the entry
  $(s,r)$.
\item  
  For $1\leq i\leq m$, if $b_{i,k}\in\{-1,1\}$, then $r_{i,k}\geq 0$.
\item  
  If there exists an entry $(2,r)$ with  $r\geq1$ in the $k^{\mbox{th}}$ column
   of the matrix $B^{ \; c,Ext}$, then at least one  entry
   $(1,r)$ must exist in the same column.  Likewise, if there exists an entry
   $(-2,r)$ with $r\geq1$
   in the $k^{\mbox{th}}$ column of the matrix
   $B^{ \; c,Ext}$, then at least one  entry $(-1,r)$ must
   exist in the same column.
\item  
 If the negative (positive) truncated  roll word of the $k^{\mbox{th}}$
 column of the matrix $B^{ \; c,Ext}$ is not an empty word, it  begins
 (ends) with an element from $\{-1,0,1\}$.
\end{enumerate}  

\item \textbf{Adjacency Properties} \par  
For $1\leq k\leq n$, the $k^{\mbox{th}}$ column of $B^{ \; c,Ext}$
satisfies these conditions.
\begin{enumerate}
\item  
  For $1\leq i\leq m$ and $2\leq k\leq n$, if $a_{i,k-1} = a_{i,k} \neq 0$,
  then $r_{i,k-1} = r_{i,k}$.
\item  
  For $1\leq i\leq m$, if $b_{i,k-1}=1$, then $b_{i,k}\in\{0,1 \}$, and
  if $b_{i,k-1}=-1$, then $b_{i,k}\in\{-1,0\}$.
  \item  
  For each ordered pair $(b,r)$ with $b \geq 2$ ($b \leq -2$) which
  appears in the $(k-1)^{\mbox{st}}$ column, there must be an index
  $i$ for which $(b_{i,k-1}, r_{i,k-1}) = (b,r)$, and $b_{i,k} \in
  C^+\cup\{1\}$ ($b_{i,k} \in C^-\cup\{-1\}$).
\item  
 For $1\leq i,j\leq m$, where $i\neq j$, if $(b_{i,k-1},r_{i,k-1}) =
 (b_{j,k-1},r_{j,k-1}) $ and  $a_{i,k}=a_{j,k}=a_{i,k-1}=a_{j,k-1}\neq 0$, then
 $b_{i,k}=b_{j,k}$.
\item  
  For $1\leq i,j\leq m$, where $i\neq j$, if $(b_{i,k-1},r_{i,k-1}) =
  (b_{j,k-1},r_{j,k-1})$, $a_{i,k}=a_{j,k}=a_{i,k-1}=a_{j,k-1}\neq 0$,
  and $b_{i,k}=b_{j,k}=b$, then there is no  $b$-factor in the word
  $b_{1,k} b_{2,k} \ldots  b_{m,k}$ which  contains both  $b_{i,k}$ and
  $b_{j,k}$.
\item  
  For every maximal 1-factor (respectively, $(-1)$-factor) $v$ in the
  word $b_{1,k} b_{2,k} \ldots  b_{m,k}$, exactly one of the following
  three conditions is fulfilled:
  \begin{enumerate} 
  \item $v$ either contains the letter $b_{m,k}$ (resp., $b_{1,k}$),
    or
  \item in the $(k-1)^{\mbox{st}}$ column there is exactly one letter
    $b_{i,k-1}=1$ (resp., $b_{i, k-1}=-1$) for which $b_{i,k }\in v$,
    or
  \item there exists exactly one sequence $v=v_1,v_2,\ldots,v_p$,
    where $p>1$, of different maximal 1-factors (respectively,
    $(-1)$-factors) in the word $b_{1,k} b_{2,k} \ldots b_{m,k}$
    satisfying the following conditions:
    \begin{itemize}
    \item For every $i$ ($1\leq i\leq p-1$) in the word
      $b_{1,k-1} b_{2,k-1} \ldots b_{m,k-1}$, there is exactly one
      letter $b_{j_i,k-1}\in\cC^+$ (resp., $b_{j_i,k-1}\in\cC^-$) for
      which $b_{j_i,k}\in v_i$, and there is exactly one letter
      $b_{s_{i+1},k-1}\in\cC^+$ (resp., $b_{s_{i+1},k-1}\in\cC^-$) for
      which $b_{s_{i+1},k}\in v_{i+1}$,
        \[ (b_{j_{i},k-1},r_{j_{i},k-1})
        = (b_{s_{i+1},k-1},r_{s_{i+1},k-1}), \]
      and $j_i \neq s_i$ for $1 < i < p$.
    \item The factor $v_p$ contains either the letter $b_{m,k}$
      (resp., $b_{1,k}$), or in the $(k-1)^{\mbox{st}}$  column there exists
      exactly one letter $b_{i, k-1}=1$ (resp., $b_{i, k-1}=-1$) for
      which $b_{i,k}\in v_p$.
    \end{itemize}
  \end{enumerate}
\item  
  For $b\notin\{-1,0,1\}$, if $v$ and $u$ represent two different
  maximal $b$-factors in the word $b_{1,k} b_{2,k} \ldots b_{m,k}$
  with the same roll number, then there is a unique sequence $ v =
  v_1,v_2,\ldots, v_p=u$, where $p>1$, of distinct maximal $b$-factors
  for which it is true that:
  \begin{itemize}
  \item For every $i$, where $1\leq i\leq p-1$, there is exactly one
    $b_{j_{i},k-1}$ with $a_{j_{i},k-1} = a_{j_{i},k} = 1$
    ($a_{j_{i},k-1} = a_{j_{i},k} = -1$), such that $b_{j_{i},k} \in
    v_{i}$, and there is a unique $b_{s_{i+1},k-1}$ with
    $a_{s_{i+1},k-1} = a_{s_{i+1},k}=1$ ($a_{s_{i+1},k-1} =
    a_{s_{i+1},k}=-1$) such that $b_{s_{i+1},k} \in v_{i+1}$,
      \[ (b_{j_{i},k-1},r_{j_{i},k-1})
      = (b_{s_{i+1},k-1},r_{s_{i+1},k-1}), \]
   and $j_i \neq s_i$.
   \end{itemize}
\end{enumerate}  

\item \textbf{Buckle Properties}  
  \quad (Specific Properties of the First, Second, and
  Last Columns)
\begin{enumerate}
\item  
  We have $b_{1,1}>0$, $b_{2,1}>0$, $(b_{1,n},r_{1,n}) =
 (b_{1,2},r_{1,2})= (0,0)$.
\item  
  For $2\leq i\leq m$, if $a_{i,n}=a_{i,1}\neq0$, then $r_{i,1} =
  r_{i,n}+1$.
\item  
  If there exists $i\in\{1,2,\ldots,m\}$ such that $r_{i,1}<0$, then
  there exists $j\in\{3,4,\ldots,m\}$ such that $j\neq i$, $a_{j,1}=
  a_{i,1}$, and $r_{j,1}=0$.
\item  
  If there exists a maximal 1-factor $b_{i_1,1} \ldots b_{i_2,1}$,
  where $1\leq i_1\leq i_2\leq m$, with $r_{i_1,1}= \ldots =
  r_{i_2,1}=0$, then
  \begin{itemize}
  \item $i_2=m$, or
  \item there exists $j_1$ with $i_2 \leq j_1 <m$  such that the word
    $b_{j_1,1}\ldots b_{m,1}$ is a maximal 1-factor with $r_{j_1,1}=
    r_{j_1+1,1}= \ldots = r_{m,1}=0$, and there exist
    $i,j\in\{1,2,\ldots,m\}$ such that $i_1\leq i\leq i_2<j_1\leq
    j<m$, $b_{i,n}= b_{j,n}$, and $r_{i,n}= r_{j,n}=-1$.
  \end{itemize}
  In addition, the first column does not contain any of the
  $-1$-factor $b_{i_1,1} \ldots b_{i_2,1}$, where $1\leq i_1\leq
  i_2\leq m$, with $r_{i_1,1}= \ldots = r_{i_2,1}=0$.
\item  
  If the last column of $B^{ \; c,Ext}$ contains the entry $(2,0)$, then
  it must  contain the entry $(1,0)$ just as well.
  Similarly, if the last column of $B^{ \; c,Ext}$ contains the entry
  $(-2,0)$, then it must  contain  the entry
  $(-1,0)$, too.
\end{enumerate}  

\item \textbf{Topological Properties}  
  \begin{enumerate}
\item  
  For $1<i_1<j_1<i_2<j_2<m$, if $b_{i_1,k}= b_{i_2,k}<-1$,
  $b_{j_1,k}=b_{j_2,k}<-1$, and $r_{i_1,k}=r_{i_2,k}=r_{j_1,k}=r_{j_2,k}$,
  then $b_{i_1,k}=b_{j_1,k}$.  Likewise, for $1\leq i_1<j_1<i_2
  <j_2 < m$, if $b_{i_1,k}=b_{i_2,k}>1$, $b_{j_1,k}=b_{j_2,k}>1$,
  and $r_{i_1,k}=r_{i_2,k}=r_{j_1,k}=r_{j_2,k}$, then $b_{i_1,k}=
  b_{j_1,k}$.
\item  
  For $1\leq i_1<j<i_2<m$, if $b_{i_1,k}= b_{i_2,k}\leq-1$,
  $b_{j,k}=-1$, and $r_{i_1,k}=r_{i_2,k}=r_{j,k}$, then $b_{i_1,k}=
  b_{i_2,k}=-1$.  Likewise, for $1\leq i_1<j<i_2\leq m$, if
  $b_{i_1,k}= b_{i_2,k}\geq1$, $b_{j,k}=1$, and $r_{i_1,k}=r_{i_2,k}
  =r_{j,k}$, then $b_{i_1,k}=b_{i_2,k}=1$.
\item  
  Assume $1\leq i<j\leq m$, if $b_{i,k}= b_{j,k}=-1$ and $r_{j,k}\neq
  r_{i,k}$, then we must have $r_{i,k} <r_{j,k}$.  Likewise, if
  $b_{i,k}=b_{j,k}=1$ and $r_{j,k}\neq r_{i,k}$, then we must have
  $r_{i,k}>r_{j,k}$.
\item  
  The absolute value of the difference between two adjacent letters in
  the negative (or positive) truncated roll word (unless it is  an
  empty word) corresponding to the $k^{\mbox{th}}$ column of the matrix
  $B^{ \; c,Ext}$ is at most~1.
\item  
  For $1\leq i,j\leq m$, if $b_{i,k}=-1$ and $b_{j,k}=1$, then $i<j$.
\item  
  If the word $b_{1,k}\ldots b_{m,k}$ does not contain 1 or $-1$,
  with the exception of eventual roots ($b_{1,k}=-1$ and/or
  $b_{m,k}=1$), and if among all the entries of the $k^{\mbox{th}}$ column of
  $B^{ \; c,Ext}$ with the same fixed roll number $r$ (note that
  $r\leq0$) there exist both negative $b_{i,k}$, where $1<i <m$,
  and positive $b_{j,k}$, where $1\leq j<m$, then the first occurrence
  of the entry $(b_{i,k},r)$ in the column with the smallest negative
  number $b_{i,k}$, such that $r_{i,k}=r$, must appear before (when
  viewed from the top row to the bottom row) the last occurrence of
  the entry $(b_{j,k},r)$ with the largest positive number $b_{j,k}$
  such that $r_{j,k}=r$.
\end{enumerate}  
\end{enumerate}  
\end{thm}

\proof
If we were to compare the statements of this Theorem, except for 5(e)
and 5(f), to the corresponding ones in Theorem $4$ of~\cite{BKDP},
which relate to HC$^{\; nc}$, we would find their formulations fairly
similar to one another. The proofs of them are thus analogous to their
 counterparts, and shall not be restated. Instead, we
move on to the two remaining exceptional cases.

Proof of 5(e): Suppose, on the contrary, that $i>j$. Then, the
shortest path in the IR from the window $w_{i,k}^0$ to its root (a
part of an NT) must cross the shortest path in the IR from the window
$w_{j,k}^0$ to its root (a part of a PT), which is impossible.

Proof of 5(f): Suppose, on the contrary, that $i>j$. Then, the
shortest path in the BRI from the positive window $w^r_{j, k}$ to its
root (located to the right and below the window $w^r_{j, k}$) must
cross the shortest path in the BRI from the negative window
$w^r_{i,k}$ to its root (located to the right and above the window
$w^r_{i,k}$), which is impossible.
\qed

\vspace*{0.5cm}

Having Part I in mind, it now comes as no surprise that Properties
1--4 are sufficient when it comes to determining a unique HC$^{ \;
  c}$.  Again, the following proof is analogous to its counterpart
from Part I. Nevertheless, in order to make this paper as
self-contained as possible, we will still provide a rough sketch of
the proof.

\begin{thm}  \label{thm:thirdEC}
Every matrix $B^{ \; c,Ext} = [(b_{i,j},r_{i,j})]_{m \times n}$ with
entries from $\big(\cC^+\cup \{1,0,-1\}\cup\cC^-\big) \times
\big\{-\big\lfloor \frac{m}{2} \big\rfloor, \ldots, \big\lfloor
\frac{m}{2}\big\rfloor\big\}$ which satisfies Properties~1--4
determines a unique HC$^{ \; c}$ on the graph $P_{m+1}\times C_{n}$.
\end{thm}

\proof The support of matrix $[b_{i,j}]_{m \times n}$, namely matrix
$[a_{i,j}]_{m \times n}$, satisfies the conditions FL$^{ \; c,Ext}$
and AC$^{ \; c,Ext}$ (Property 1(a) of Theorem~\ref{thm:secondC}).  We
will prove that RC$^{ \; c,Ext}$ holds, through a set of claims.  But
first, take the set of all windows of $\cG_m$ into consideration and
divide them into positive, negative and zero ones, in accordance with
the sign of the value corresponding to $b_{i,j}$.  Note that the
window corresponding to $1$ can not be adjacent to a window
corresponding to $-1$ because of Property 1(a).

The edges of $\cG_m$ which belong to different kinds of windows (zero
and non-zero ones), together with the edges of the zero windows
belonging to lines $M_0M_1$ and $N_0N_1$, determine a spanning
$2$-regular subgraph of $\cG_m$.  Adding the lines $M_0M_1$ and
$N_0N_1$ to it gives way to a clear distinction between the positive,
negative and zero regions.  The first, of course, being determined by
$b_{i,j} > 0$, the second by $b_{i,j} < 0$, and the last one by
$b_{i,j} = 0$.  However, instead of focusing on these regions per say,
we can observe the components of the subgraph $\cW_m$ induced by the
windows of the same kind (positive, negative or zero).  We will refer
to them as the \textbf{\emph{ positive, negative}} or \textbf{\emph{
zero}} regions ``induced by the positive, negative or zero entries
of the matrix $[b_{i,j}]_{m \times n}$''.

Note that every entry $(b_{ij},r_{ij})$ of the matrix $B^{ \; c,Ext}$
is assigned to \emph{exactly one} window $w_{i,j}^r$, where $r =
r_{ij}$. In other words, $w_{i,j}^r$ belongs to the
$(j+nr)^{\mbox{th}}$ column of $\cG_m$ (the square
$M_rN_rN_{r+1}N_{r+2}$), although there are infinitely (countably)
many vertices of $\cW_m$ corresponding to this $b_{ij}$.  If we
collect all the positive and negative windows assigned to entries of
the matrix $B^{ \; c,Ext}$ we will obtain a finite number of
completely fulfilled regions, as the claim below shows.

\vspace*{0.2cm}

\begin{clm}  \label{clm:first}
Every window from any positive or negative region that contains
$w_{i,j}^r$, where $r = r_{ij}$, is assigned to an entry of the matrix
$B^{ \; c,Ext}$.
\end{clm}

\noindent\emph{Proof.}
Since there is a path between any two windows in the considered
regions, this  comes as a consequence of Properties 1(b), 3(a) and
4(b).
\qed

\vspace*{0.2cm}

As a result, every positive or negative region is bounded, and there
are infinitely (countably) many regions congruent to it.  The regions
described in the previous lemma will be called the \textbf{\emph{basis
positive regions}} or the \textbf{\emph{basis negative regions}}
and its each of its windows $w_{ij}^r$'s for which $(b_{ij},r) =
(b,r)$ a \textbf{\emph{$b$-window}}, where $b \neq 0$.  If any such
window belongs to the last (that is, the $m^{\mbox{th}}$) row and $b
\neq 0$, it must be a $1$-window with $r=0$ (Property 1(d)); it will
be called the \textbf{\emph{ down root}}.  If any such window belongs
to the first row and $b<0$, it must be a $-1$-window with $r=0$
(Property 1(d)); it will be called the \textbf{\emph{ up root}}. The
window $w_{11}^r$, where $r = r_{11}$ is the only \textbf{\emph{ up
root}} which is a positive window.  Recall that in case the path
which connects the window $w_{i,j}^r$ to another window consists only
of windows from its column (the $(j+nr)^{\mbox{th}}$ column) or/and
those to the left of it, we call this path the \textbf{\emph{left path
for the window $w_{i,j}^r$}}.

\vspace*{0.2cm}

\begin{clm}  \label{clm:second}
For any two windows $w_{i,j}^r$ and $w_{i',j}^r$, where $i <i'$, from
the same basis positive (or negative) region and the same column (the
$(j+nr)^{\mbox{th}}$ column) for which there exists a left path for
and between them, the following must be fulfilled: $b_{i,j} =
b_{i',j}$.
\end{clm}

\noindent\emph{Proof.}
This can be proved by strong induction on the length $l$ of the
considered path using Properties~1(b) and~3(d), in the exact same way
we did in the proof of Lemma 2 in \cite{BKDP}.
\qed

\vspace*{0.2cm}

\begin{clm}  \label{clm:third}
The subgraph of $W_{m,n}$ induced by positive (or negative) entries of
matrix $B^{c,Ext}$ has a forest  structure.
\end{clm}

\noindent\emph{Proof.}
Assuming the opposite holds, that there exists a cycle in a basis
positive region, then in the rightmost column of its windows once we
apply Claim~\ref{clm:second} we reach a contradiction with either
Property 3(e), 3(f) or 3(g) (compare with the proof of Lemma 3 in
\cite{BKDP}).
\qed

\begin{clm}  \label{clm:fourth}
Let $w_{i_1,j}^r$ and $w_{i_2,j}^r$ be any two windows from the same
basis positive region, which belong to the same column (the
$(j+nr)^{\mbox{th}}$ column), with $b_{i_1,j} = b_{i_2,j}=b$ and $\mid
b \mid >1$ ($r_{i_1,j} = r_{i_2,j}=r$). Then, there exists a unique
left path for and between them in this region.
\end{clm}

\noindent\emph{Proof.}
The existence of such a left path is proved by induction on $j+nr$
using Property 3(g) and ~3(b).  The base case deals with the leftmost
windows of the considered region, whereas Claim~\ref{clm:third}
implies its uniqueness.
\qed

\vspace*{0.2cm}

\begin{clm}  \label{clm:fifth}
For every $1$-window (resp., $(-1)$-window) $w_{i,j}^r$, where $r =
r_{ij}$, there exists a unique left path for it which connects it to a
down root (resp., an up root).
\end{clm}

\noindent
\emph{Proof.}
The proof can be obtained by induction on $j+nr$.  If the window
$w_{i,j}^r$ belongs to the leftmost windows in the considered region
(the base case of the induction), the letter $b_{i,j}$ and $b_{m,j}$
(resp., $b_{1,j}$) must belong to the same $1$-factor (resp.,
$(-1)$-factor) (Property 3(f)i).  If it is not the case, from Property
3(f) and Claim~\ref{clm:fourth} we conclude that either there is a
unique left path for and from it to $w_{m,j}^r$ ($w_{1,j}^r$) which is
a down root (resp., an up root), or there is a unique left path for
and from it to a unique $1$-window (resp., $(-1)$-window) from the
previous column (the $(j+nr-1)^{\mbox{th}}$ column).  In the second
case, we apply the induction hypothesis to the newly obtained
$1$-window (resp., $(-1)$-window) instead of to $w_{i,j}^r$.
\qed

\vspace*{0.2cm}

\begin{clm}  \label{clm:sixth}
Every positive region has a unique down root, whereas every negative
region has a unique  up root.
\end{clm}

\noindent
\emph{Proof.}
Property 3(c)  implies that every  rightmost window  of any basis
positive (resp., negative) region is a $1$-window (resp.,
$(-1)$-window).  By applying Claim~\ref{clm:fifth} to these windows we
obtain the desired statement.
\qed

\vspace*{0.5cm}

Now we can finish the proof of the main statement.
Claims~\ref{clm:sixth} and~\ref{clm:third} together with Property 3(c)
imply that the RC$^{c,Ext}$ is satisfied. By applying
Lemma~\ref{lem:firstC} we finally obtain the existence and uniqueness
of a HC$^{\; c}$ on the graph $P_{m+1}\times C_{n}$ whose split tree
contains the window $w_{1,1}$.
\qed

\vspace*{0.5cm}

For each integer $m\geq1$, we will create an auxiliary digraph whose
role will be to enumerate the number of HC$^{\; c}$'s in $P_{m+1}
\times C_n$.  Here is how we intend to do that.  At first, let $\cF_m
= \cF^{ \; c,Ext}_m$ denote the set of all the possible first columns
of $B^{c,Ext}$, and $\cD^{c,Ext}_m$ a digraph with the vertex set
$V({\cal D}^{c,Ext}_{m})$ which consists of all the possible remaining
columns of the same matrix. For any $v,u \in V({\cal D}^{c,Ext}_m)$,
there exists an arc from $v$ to $u$ if and only if the vertex
  \[ v = [(b_{1,k},r_{1,k}), (b_{2,k},r_{2,k}),
  \ldots, (b_{m,k},r_{m, k})]^T \]
may appear as a column preceding the vertex
  \[ u = [(b_{1,k+1},r_{1,k+1}), (b_{2,k+1},r_{2,k+1}),
  \ldots, (b_{m,k+1}, r_{m,k+1})]^T, \]
for $2 \leq k \leq n-1$.  Note that the vertices of the disjoint sets
$\cF_m$ and $V({\cal D}^{c,Ext}_{m})$ are in both cases the column
vectors of the form $[(b_1,r_1), (b_2,r_2), \ldots , (b_m,r_m)]^T$
with entries from $(\cC^-\cup \{-1,0,1\}\cup \cC^+ ) \times \big\{
-\lfloor\frac{m}{2}\big\rfloor,\ldots,\big\lfloor\frac{m}{2}\big
\rfloor\big\}$.  The difference between the two is that $ b_1 = b_2 >
0$ for the vertices from $\cF_m$, whereas $b_1 \in \{ -1, 0 \} $ and
$r_1 = 0$ for the vertices in $V({\cal D}^{c,Ext}_{m})$.

Let $\cS_m, \cL_m \subseteq V({\cal D}^{c,Ext}_{m})$ denote the set of
all possible second and last (that is, the $n^{\mbox{th}}$) columns of
the matrix $B^{c,Ext}$, respectively.  Also, let $\cL \cF \cS_m$
denote the set of all possible ordered triples $(l,f,s) \in \cL_m
\times \cF_m \times \cS_m$ of columns which can appear as the last
($n^{\mbox{th}}$), first and second column, respectively, in
$B^{c,Ext}$.  The aforementioned auxiliary diagraph from the previous
paragraph will be denoted by $ \overline{{\cal D}}^{c,Ext}_{m}$. Its
set of vertices will be $V(\overline{{\cal D}}^{c,Ext}_{m})=\cF_m \cup
V({\cal D}^{c,Ext}_{m})$ and its set of edges
  \begin{eqnarray*}
  \lefteqn{E(\overline{{\cal D}}^{c,Ext}_{m})} \\
  &=& E({\cal D}^{c,Ext}_{m}) \cup \{ (u,v) \mid
  (\exists w) (u,v,w) \in \cL \cF \cS_m \vee
  (\exists w) (w,u,v) \in \cL \cF \cS_m \}.
  \end{eqnarray*}
Note that all the vertices of this graph do abide by the Basic and
Column properties (as well as by the Topological properties).
Additionally, the arcs of the digraph ${\cal D}^{c,Ext}_{m}$ abide by
the Adjacency properties, whereas the arcs coming out of the set
$\cL_m$ and into the vertices from the set $ \cF_m$ satisfy both the
Adjacency and Buckle properties. The same goes for the arcs spanning
from the set $ \cF m$ and into the vertices from the set $ \cS_m$

\begin{figure}[htb]
\begin{center}
\includegraphics{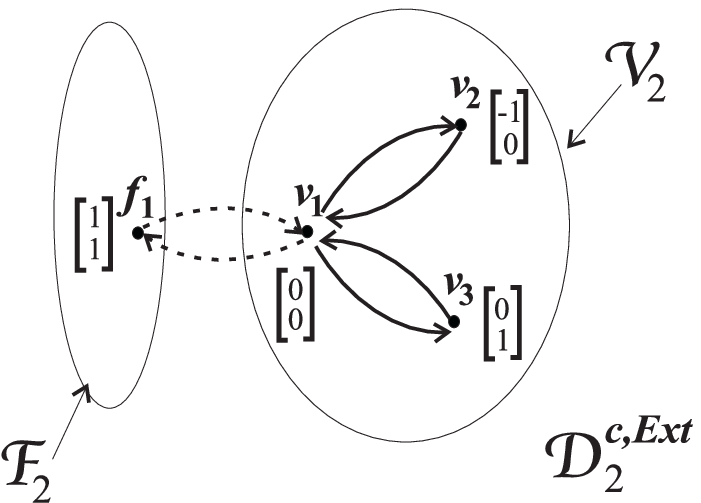}
\end{center}
\caption{The digraph $\cD^{ \; c,Ext}_2$ and the corresponding set
  $\cF_2$.}
\label{fig:caseEm=2}
\end{figure}

For example, when $m=2$, the digraph $ \overline{{\cal D}}^{c,Ext}_{2}$
has four vertices, and the set $\cL \cF \cS_m$ consists of just one
triplet $(v_1,f_1,v_1)$ (see Figure~\ref{fig:caseEm=2}). When $m=3$,
the digraph $ \overline{{\cal D}}^{c,Ext}_{3}$ has fourteen vertices,
3 of which belong to the set $\cF_3 = \{ f_1, f_2, f_3 \}$; whereas 11
of them as in $V({\cal D}^{c,Ext}_{3})$. At the same time, there exist
precisely two arcs from each of the vertices from the set $\cF_3 $
into the set $\cS_3 = \{ v_1, v_2, v_3 , v_4, v_5 \}$. Also, there
exist two arcs per every vertex of the set $\cF_3$ to which they point
from the set $\cL_3 = \{ v_1, v_2, v_3 , v_{10} \}$ thus forming 12
triplets - elements of the set $\cL \cF \cS_3 $ (see Subsection 5.3).

In this way, the enumeration of HC$^{\; c}$s on $P_{m+1} \times C_n$
is reduced to the enumeration of oriented walks of length $n-2$ in the
digraph $\cD^{c,Ext}_m$ with the pairs of initial and final vertices
which are respectively the third and first coordinates of the triplets
from the set $\cL \cF \cS_m $. In other words, this enumeration is
reduced to the enumeration of closed oriented walks of length $n$ in
the digraph $\overline{\cD}^{c,Ext}_m$ for which it holds that they
both start and finish in the same vertex from the set $\cF_m$ and no
other vertex from the set $\cF_m$ belongs to them.  Finally, this
number $\varphi^{c, Ext}_m(n-2)$, where $n\geq 2$, needs to be
multiplied by $n$ so as to obtain the correct number of HC$^{ \; c}$
of $P_{m+1} \times C_n$.


\vspace{1.0cc}
\begin{center}
{\bf 3. CODING THE INTERIOR TREE BY NON-ZERO ENTRIES}
\end{center}

\vspace{0.5cc}
\begin{center}
{\bf 3.1.  The First Phase --- the Matrix $A^{ \; c,Int}$}
\end{center}

Here, the zero windows belong to the exterior trees and $w_{11}$
remains the up root of the split tree.
To put it differently,
$w_{1,2}^0$ is the leftmost window from the first row of the interior
region in the BRI.

Each HC$^{\; c}$ on $P_{m+1}\times C_n$, where $m\geq1$, with the
window $w_{1,1}$ as the up root of the split tree can be encoded by a
$(0,1)$-matrix $A^{ \; c,Int} = [a_{i,j}]_{m \times n}$ where
  \[ a_{i,j} \stackrel{\rm def}{=} \begin{cases}
  1 &  \mbox{ if $w_{i,j}$ belongs to the interior of HC$^{\; c}$}, \cr
  0 &  \mbox{ otherwise.}  \cr
  \end{cases} \]
By doing so, we obtain one positive region and one or more zero
regions in the BRI.  Note that in Figure~\ref{fig:int-flat} almost all
the windows of the IT are in ${\cal R}_0$ except for the two windows
which are in ${\cal R}_1$.  On the other hand, in Figure~\ref{fig:RI}
the windows of the IT belong to the rectangles ${\cal R}_{-2}, {\cal
R}_{-1}$ and $ {\cal R}_0$.

\begin{lem}  \label{lem:firstIC}
Every HC$^{\; c}$ on the thick grid cylinder graph $P_{m+1} \times
C_n$ with the window $w_{1,1}$ as the up root of the split tree
determines a $(0,1)$-matrix $A^{ \; c,Int} = [a_{ij}]_{m \times n}$
that satisfies the following conditions ($a_{i,n+1} \stackrel{\rm
def}{=} a_{i,1}$, and $a_{i,0} \stackrel{\rm def}{=} a_{i,n}$ for
$1\leq i\leq m$).

\newpage

\begin{enumerate}
\item \textbf{First and Last Row Conditions} (FL$^{ \; c,Int}$):
  \begin{enumerate}
  \item $a_{1,1}=0$.
  \item For $1\leq   j\leq n$,  $(a_{1,j},a_{1,j+1})    \not\equiv(0,0)$.
  \item For $1\leq  j\leq n$,  $(a_{m,j},a_{m,j+1})    \not\equiv(0,0)$.
  \end{enumerate}

\item \textbf{Adjacency Conditions} (AC$^{ \; c,Int}$):
For $1\leq  i \leq m-1$, and $1\leq  j \leq n$,
  \[ (a_{i,j},a_{i+1,j},a_{i,j+1},a_{i+1,j+1})
  \not\in \{(1,1,1,1),(0,0,0,0),(0,1,1,0),(1,0,0,1)\}. \]

\item \textbf{Tree Condition} (TC$^{ \; c,Int}$):
The vertices of  $W_{m,n}$ corresponding to 1's in $A^{ \; c,Int}$ induce a
unique tree in $W_{m,n}$.

\end{enumerate}
Conversely, every $(0,1)$-matrix $[a_{ij}]_{m \times n}$ which
satisfies the  conditions FL$^{ \; c,Int}$, AC$^{ \; c,Int}$, and
TC$^{\; c,Int}$ determines a unique HC$^{\; c}$ on the thick grid
cylinder graph $P_{m+1} \times  C_n$  with the window $w_{1,1}$ as the
up root of its split tree.
\end{lem}

\proof
The first two conditions provide the local whereas the third one
provides the global aspect of hamiltonicity and their necessity is
easily verifiable (note that FL$^{ \; c,Int}$ implies that
$a_{1,2}=a_{1,n}=1$).  With the intention of showing that all the
above mentioned conditions are sufficient as well, note the following.
The first two conditions ensure that the set of edges belonging to
both a zero and a positive window or to both a positive window and one
of the lines $M_0M_1$ or $N_0N_1$ determines a unique $2$-regular
spanning subgraph of $P_{m+1} \times C_n$, that is, a union of
cycles. The third condition implies that there exists a unique cycle
--- the boundary of the positive region (IT).

Let us walk from the upper horizontal edge of the window $w_{12}$ (the
windows on the right-hand side belong to the IT), walking in
accordance with the aforementioned boundary.  As there are no
consecutive zeros in the $m^{\mbox{th}}$ row of $A^{ \; c,Int}$, there
exists at least one window corresponded to $1$ in that row.  Find the
last lower horizontal edge of some window from the $m^{\mbox{th}}$ row
through which we pass along our walk.  If we denote that window by
$w_{mk}$ (in Figure~\ref{fig:int-flat}, we have $k=4$), then
$w_{m,k-1}$ represents the down root of the split tree. The reason
behind that is that the rest of our walk consists of edges which
belong to zero windows that are connected to both $w_{m,k-1}$ and
$w_{11}$, with the latter of which we end our walk as it is.
Therefore, the constructed HC is contractible and has $w_{1,1}$ for
the up root of the split tree.  \qed

\vspace{0.5cc}
\begin{center}
{\bf 3.2. The Second Characterization of HC$^{\; c}$ with $w_{11}$ \\
as  the Up Root of the Split Tree}
\end{center}

For each HC$^{\; c}$ with $w_{1,1}$ in the split tree, we associate
the matrix $A^{ \; c,Int} = [a_{i,j}]_{m \times n}$ to the matrix $B^{
\; c,Int} = [(b_{i,j},r_{i,j})]_{m\times n}$.
The first matrix satisfies the conditions FL$^{ \; c,Int}$, AC$^{ \;
  c,Int}$, and TC$^{ \; c,Int}$.  The second matrix with $b_{i,j}
\in\cC^+ \cup\{0,1\}$ and $- \big\lfloor \frac{m}{2} \big\rfloor \leq
r_{i,j} \leq \big\lfloor \frac{m}{2} \big\rfloor$ is constructed in
the following way:

\begin{enumerate}
\item Define $r_{i,j} = r(w_{i,j})$.
\item Set  $b_{i,j} = a_{i,j}=0$ if $w_{i,j}$ belongs to an ET.
\item For each fixed column $j$, partition the  positive windows from
  the $j^{\mbox{th}}$ column with the same roll number into
  $j^r$-joined equivalence classes.  Then, label all the windows within
  each equivalence class with $2$, $3 \ldots $, according to the order
  in which the equivalence classes first appear within the
  $j^{\mbox{th}}$ column, from top to bottom.
\end{enumerate}

\vspace{0.5cc}
\noindent\emph{Example 6.} \ %
In Figure~\ref{fig:int-flat}, the values $(b_{ij},r_{ij})$ of $B^{ \;
c,Int}$ are inscribed on the windows as $\ds b_{ij}^{r_{ij}}$ or just
as $\ds b_{ij}$ if $r_{ij}=0$.  In the $9^{\mbox{th}}$ column
there are four parts of the IT that belong to the roll~0.  Three of
them, $w_{1,9}$, $w_{4,9}$, and $w_{8,9}$, are $9^{0}$-joined.
Consequently, the same $b$-value is assigned to them in $B^{ \;
c,Int}$.  More specifically, we have $b_{1,9}=b_{4,9}=b_{8,9}=2$.
The fourth window, $w_{6,9}$, while still belonging to the roll~0,
belongs to a different equivalent class.  Hence, $b_{6,9}=3$.  In the
second column, there are three windows that belong to the IT.  Two of
them belong to the roll~0, but they belong to two different
equivalence classes with respect to the relation $2^{0}$-joined.  The
last window, $w_{10,2}$, belongs to roll~1; thus, $(b_{10,2},r_{10,2})
= (2,1)$.

\vspace{0.5cc}

In an arbitrary column of the matrix $[b_{i,j}]_{m \times n}$ that
corresponds to $B^{ \; c,Ext} = [(b_{i,j},r_{i,j})]_{m\times n}$, we
consider all the maximal $b$-factors (if it exists), where $b>0$.  Let
them be, in their order of appearance, that is, from top to bottom,
$p_1$-factor, $p_2$-factor, \ldots, $p_k$-factor, where $k\geq1$, and
$p_i\geq2$ for each $i$.  In addition, let $r_1,r_2,\ldots,r_k$ denote
the roll numbers associated with these maximal factors. The words $p_1
p_2 \ldots p_k$ and $r_1 r_2 \ldots r_k$ are called the \textbf{\emph{
truncated word}} and the \textbf{\emph{truncated roll word}},
respectively. A subsequence of a truncated word induced by the letters
with the same roll number $r$ is called a \textbf{\emph{color$^r$
word}}.

\vspace{0.5cc}
\noindent\emph{Example 7.} \ %
For the second column in Figure~\ref{fig:int-flat}, the truncated
word, the truncated roll word, the color$^0$ word, and the color$^1$
word are: 232, 001, 23, and 2, respectively.  The color$^0$ word for
the fourth column is 23435.

\vspace{0.5cc}
\begin{center}
{\bf 3.3. Properties of $B^{ \; c,Int}$}
\end{center}

The properties listed below follow straightforwardly from the
definition of the matrix $B^{ \; c,Int} = [(b_{i,j},r_{i,j})]_{m
\times n}$.  Here, $(b_{i,n+1},r_{i,n+1}) \stackrel{\rm def}{=}
(b_{i,1},r_{i,1})$, and $(b_{i,0},r_{i,0}) \stackrel{\rm def}{=}
(b_{i,n}, r_{i,n})$.

\begin{thm}  \label{thm:secondIC}

The matrix $B^{ \; c,Int} = [(b_{i,j},r_{i,j})]_{m \times n}$
satisfies the following conditions.
\begin{enumerate}

\item \textbf{Basic Properties}  
\begin{enumerate}
\item  
  The support of the matrix $[b_{i,j}]_{m \times n}$, that is, the
  matrix $[a_{i,j}]_{m \times n}$, satisfies the
  FL$^{ \; c,Int}$ and AC$^{ \; c,Int}$.
\item  
  \emph{Harmonization of the adjacent entries having the same sign}:
  For $2\leq i\leq m$ and $1\leq k\leq n$, if $a_{i-1,k}=a_{i,k}$,
  then $(b_{i-1,k},r_{i-1,k}) = (b_{i,k},r_{i,k})$.
 \item  
  For $1\leq i\leq m$ and $1\leq j\leq n$, if $a_{i,j}=0$, then
  $r_{i,j}=0$.
\end{enumerate}  

\item \textbf{Column Properties} \par
  For $1\leq k\leq n$, the $k^{\mbox{th}}$ column
  $[(b_{1,k},r_{1,k}),(b_{2,k},r_{2,k}),\ldots,(b_{m,k},r_{m,k})]^T$
  of the matrix $B^{ \; c,Int}$ satisfies these conditions:
\begin{enumerate}
\item  
  If there exists an entry $(s,r)$ in the $k^{\mbox{th}}$ column of
  $B^{ \; c,Int}$, where  $s\geq3$, then for each $\ell\in\{2,
  3,\ldots,s-1\}$, at least one copy of the entry $(\ell,r)$ must
  appear before the first appearance of the entry $(s,r)$.
\item  
  If the  truncated roll word of the first column of $B^{ \; c,Int}$ is
  not an empty word, it begins with $0$ or $1$. The  truncated roll
  word of the $k^{\mbox{th}}$ column of $B^{ \; c,Int}$  for  $k \neq
  1$ is non-empty  and begins with $0$.
\end{enumerate}  

\item \textbf{Adjacency Properties}  

For $1\leq k\leq n$, the $k^{\mbox{th}}$ column of $B^{ \; c,Int}=
[(b_{i,j},r_{i,j})]_{m \times n}$ satisfies these conditions:
\begin{enumerate}
\item  
  For $1\leq i\leq m$ and $2\leq k\leq n$, if $a_{i,k-1}=a_{i,k}=1$,
  then $r_{i,k-1} = r_{i,k}$.
\item  
  For each ordered pair $(b,r)$ that appears in the $k^{\mbox{th}}$
  column, with the exception of the following two cases:
  \begin{itemize}
  \item $b=2$, $r>0$ is the maximal roll in this column, and there
    is no occurrence of $(3,r)$ in this column;
  \item $k=n$, $b=2$, $r=0$ is the maximal roll in this column, and
    there is no occurrence of $(3,0)$ in this column;
  \end{itemize}
  there must exist an index $i$ for which $(b_{i,k},r_{i,k})=(b,r)$,
  and $b_{i,k+1}\neq0$.

\item  
  For $1\leq i,j\leq m$, where $i\neq j$, if $(b_{i,k-1},r_{i,k-1}) =
 (b_{j,k-1},r_{j,k-1}) $ and  $a_{i,k}=a_{j,k}=a_{i,k-1}=a_{j,k-1} =1$, then
 $b_{i,k}=b_{j,k}$.
\item  
  For $1\leq i,j\leq m$, where $i\neq j$, if $(b_{i,k-1},r_{i,k-1}) =
  (b_{j,k-1}, r_{j,k-1})$, $a_{i,k}=a_{j,k}=a_{i,k-1}=a_{j,k-1}=1$,
  and $b_{i,k}=b_{j,k}=b$, then there is no $b$-factor in the word
  $b_{1,k} b_{2,k} \ldots b_{m,k}$ which contains both $b_{i,k}$ and
  $b_{j,k}$.
\item  
  For $b \neq 0$, if $v$ and $u$ represent two different maximal
  $b$-factors in the word $b_{1,k} b_{2,k} \ldots  b_{m,k}$ with the
  same roll number, then there is a unique sequence
  $v=v_1,v_2,\ldots,v_p=u$  $(p>1)$ of  distinct maximal $b$-factors
  for which it is true that:
  \begin{itemize}
  \item For every $i$, where $1\leq i\leq p-1$, there is exactly one
    $b_{j_i,k-1}$ with $a_{j_i,k-1}=a_{j_i,k}=1$, such that  $b_{j_i,k}
    \in v_i$, and there is a unique $b_{s_{i+1},k-1}$ with
    $a_{s_{i+1},k-1}=a_{s_{i+1},k}=1$ such that $b_{s_{i+1},k}\in
    v_{i+1}$, $(b_{j_i,k-1},r_{j_i,k-1}) = (b_{s_{i+1},k-1},
    r_{s_{i+1},k-1})$, and  $j_i \neq s_i$.
  \end{itemize}

\item  
  For $r>0$, if the ordered pair $(2,r)$ appears in the $k^{\mbox{th}}$ column,
  then there must exist an index $i$ for which
  \begin{itemize}
  \item If $k=1$, then $r_{i,n}=r-1$, and $a_{i,1} =a_{i,n} = 1$.
  \item If $k>1$, then $r_{i,k-1}=r$, and $a_{i,k-1} =  a_{i,k}= 1$.
  \end{itemize}
\end{enumerate}  

\item \textbf{Buckle Properties}
  \quad (Specific Properties of the First and Last Columns) 
  \begin{enumerate}
  \item  
    $(b_{1,1},r_{1,1})=(b_{2,1},r_{2,1})=(0,0)$, and $(b_{1,n},r_{1,n})
    =(2,0)$.
  \item  
    For $2<i \leq m $, if $a_{i,n}=a_{i,1}=1$, then $r_{i,1} =  r_{i,n}+1$.
  \end{enumerate}  

\item \textbf{Topological Properties}  
  \begin{enumerate}
  \item  
  For $1<i_1<j_1<i_2<j_2\leq m$, if $b_{i_1,k}= b_{i_2,k}>1$,
  $b_{j_1,k}=b_{j_2,k}>1$, and $r_{i_1,k}=r_{i_2,k}=r_{j_1,k}
  =r_{j_2,k}$, then $b_{i_1,k} = b_{j_1,k}$.
\item  
  The absolute value of the difference between two adjacent letters in
  the non-empty truncated roll word
  corresponding to the $k^{\mbox{th}}$ column of $B^{ \; c,Int}$ is at most~1.
\end{enumerate}  
\end{enumerate}  
\end{thm}

\proof
We shall omit the proofs of those items that we consider to be fairly
straightforward, due to their similarity to the ones in
Theorem~\ref{thm:secondC} or Theorem $4$ of~\cite{BKDP}.  However, we
shall discuss the remaining items.

{\em Column Property 2(a)}:
It is a trivial consequence of the chosen method of coding, that is, of
the way in which $B^{ \; c,Int}$ is formed.

{\em Column Property 2(b)}:
Notice at first that using the definition of the matrix $B^{ \;
c,Int}$ and the definition of BRI, we have that $b_{1,1}=
b_{2,1}=0$, $b_{1,2}=b_{1,n}=2$ and $r_{1,2} = r_{1,n} = 0$.  The
number $r_{1,n}$ must be zero because, assuming the contrary, the path
which connects the windows $w_{1,2}^0$ and $w_{1,n}^{r_{1,n}}$ would
encompass the split tree and so it could not reach its down root.

The truncated roll word of the $k^{\mbox{th}}$ column where $k>1$ is a
non-empty word because $b_{1,k}b_{2,k} \ldots b_{m,k}$ can be zero
word just for $k=1$.  Now, suppose to the contrary, that the truncated
roll word of the $k^{\mbox{th}}$ column, where $k>1$, of $B^{\; c,
Int}$ begins with $r \geq 1$.  Let $w^{\; r}_{i,k}$ be a positive
window from the BRI, for which $b_{lk} = 0$, where $1 \leq l \leq
i-1$.  Let $P^0$ denote the unique path from $w^{\; r}_{i,k}$ to
$w^{\; 0}_{1,2}$.  This path must cross the $(k+n(r-1))^{\mbox{th}}$
column of the RI in a window $w^{\; r-1}_{s,k}$ for some integer
$s>i$.  But, in that case, the copy $P^{-1}$ of the path $P^0$ in the
RI meets $P^0$, which is impossible.  This proves that $r<1$.
Similarly, we conclude that $r >-1$, observing the path from $w^{\;
r}_{i,k}$ to $w^{\; 0}_{1,n}$ instead of the path from $w^{\;
r}_{i,k}$ to $w^{\; 0}_{1,2}$.  Together, we conclude that $r=0$.

Further we show that the first column of $B^{ \; c,Int}$ possesses the
property that a non-empty truncated roll word begins with $r_{i,1}=0$
or $r_{i,1}=1$ ($3 \leq i \leq m$).  Namely, assuming the opposite, we
would have $r_{i,1} >1$ or $r_{i,1} <0$. If $r_{i,1} >1$, then the
unique path $P^0$ from $w^{\; 0}_{1,n}$ to $w^{\; r_{i,1}}_{i,1}$ must
cross the $(n+1)^{\mbox{st}}$ column of the RI in a window $w^{\;
1}_{s,1}$ for some integer $s>i$. It implies that the copy $P^{-1}$
of the path $P^0$ in the RI meets $P^0$, which is impossible.
Similarly, we may consider the case when $r_{i,1} <0$ and reach a
contradiction again.

{\em Adjacency Properties 3(b)}:
Assume the rightmost positive windows in the BRI correspond to
column $\tilde{k}$.  They have maximal roll in this column, say
$\tilde{r}$, which is either (i) $\tilde{r}>0$, or (ii) $\tilde{r}=0$
and $\tilde{k}=n$. Since there exists exactly one
$\tilde{k}^{\tilde{r}}$-joined equivalence class for these windows,
all corresponding entries of  the matrix $B^{\; c, Int}$ are
$(2,\tilde{r})$.  For any other window $w_{j,k}^r$ in the BRI which
corresponds to the pair $(b,r)$ there is a path from it to one of these
rightmost positive windows in the BRI over some window $w^{r}_{i,k+1}$
from the next, $(rn+(k+1))^{\mbox{th}}$ column. Clearly,
 $(b_{i,k}, r_{i,k}) = (b,r)$.

Additionally, note that if the pair $(b_{i,j},r_{i,j}) = (2,r)$, where
$ r > 0$ is the maximal roll in the $j^{\mbox{th}}$ column of $B^{ \;
  c,Int}$, and there is no occurrence of $(3,r)$ in this column, then
the window $w_{i,j}^r$ can be (but need not be) one of the rightmost
windows in the considered region (see, as an example, the first and
the second column of $B^{ \; c,Int}$ in Figure~\ref{fig:int-flat}).

{\em Adjacency Properties 3(f)}:
If the ordered pair $(2,r)$ with $r>0$ appears in the $k^{\mbox{th}}$
column, then the corresponding window belongs to the $(rn +k)^{th}$ column
in BRI. Since there exists a path from it to the window $w_{1,2}$, it
must pass through the previous column.
\qed

\medskip
Properties 1--4 are sufficient for determining a unique HC$^{ \; c}$.

\begin{thm}  \label{thm:thirdIC}
Every matrix $B^{ \; c,Int} = [(b_{i,j},r_{i,j})]_{m \times n}$ with
entries from $\big(\cC^+\cup \{1,0\}\big) \times
\big\{-\big\lfloor \frac{m}{2} \big\rfloor, \ldots, \big\lfloor
\frac{m}{2}\big\rfloor\big\}$ which satisfies Properties~1--4  of
Theorem~\ref{thm:secondIC} determines a unique HC$^{ \; c}$ on the
graph $P_{m+1}\times C_{n}$.
\end{thm}

\proof
By applying Property 1(a), we see that the support of the matrix
$[b_{i,j}]_{m \times n}$, in other words, the matrix $[a_{i,j}]_{m
  \times n}$, satisfies the FL$^{ \; c,Int}$ and AC$^{ \; c,Int}$ of
Lemma~\ref{lem:firstIC}.  In order to prove that the TC$^{ \; c,Int}$
of Lemma~\ref{lem:firstIC} is satisfied as well, observe the set of
all non-zero windows $w_{i,j}^{k}$ ($ k \in \mathbb{Z}$) of $\cW_m$
--- the windows corresponding to the entries $(b_{i,j},r_{i,j})$ of
$B^{ \; c,Int}$ where $b_{i,j}>0$. These windows determine the {\em
non-zero regions}. The union of the boundaries of these regions
present a spanning 2-regular subgraph of $\cG_m$.  If we consider only
the non-zero windows $w_{i,j}^{k}$ where $k= r_{i,j}$, applying
Properties 1(b), 3(a) and 4(b), we can conclude that they cover up
completely one or more of those regions.  To distinguish these regions
from their copies we shall call them the \textbf{\emph{basis regions
(BR)} } (although there is just one such region, as will be shown
later).

Recall that in case a path in $\cW_m$, which connects the window
$w_{i,j}^r$ to another window, consists only of the windows from its
column (the $(j +nr)^{\mbox{th}}$ column) or/and those to the left of
it, we call it the \textbf{\emph{left path for the window
$w_{i,j}^r$}}.  The path visiting only the windows from the same
column which are assigned to the entry $(b,r)$ of $B^{c,Int}$ is
called a \textbf{\emph{\boldmath{$b$}-factor}}.  Properties~1(b)
and~3(c) imply the next statement.

\begin{clm}  \label{clm:seventh}
If  there exists a left path for and between the  two windows
$w_{i,j}^r$ and $w_{i',j}^r$, where $i<i'$, from the same basis
region, then $b_{i,j} = b_{i',j}$.
\end{clm}

The proof of this statement can be obtained by strong induction on the
length $l$ of the considered path in a similar fashion by which
Claim~\ref{clm:second} was proven. Further, using
Claim~\ref{clm:seventh}, Properties~3(d) and~3(e), analogously as in
the proof of Claim~\ref{clm:third}, the following statement can be
shown.

\begin{clm}  \label{clm:eighth}
The subgraph of $W_{m,n}$ induced by the windows determined by the
positive  entries of $B^{c,Int}$ has  a forest  structure.
\end{clm}


The next claim can be proved by induction on $j+nr$ (with the base
case which refers to the leftmost windows of the considered region),
whilst at the same time relying upon Property 3(e).

\begin{clm}  \label{clm:nineth}
For any two windows $w_{i_1,j}^r$ and $w_{i_2,j}^r$  from the same
basis  region and with  $b_{i_1,j} = b_{i_2,j}$  there  exists a left
path for and between them.
\end{clm}

The uniqueness of this path is a consequence of the forest structure
of the subgraph of $W_{m}$ induced by the windows belonging to the
considered region.


It remains to prove that the subgraph of $W_{m}$ induced by the
positive windows which belong to the basis regions has just one
component (hence it is a connected graph).  We will prove that there
exists a path for and between an arbitrary such window and
$w_{1,n}^{0}$ (Property 4(a) guaranties the existence of such a path
in the considered region).

Note that Property 3(b) implies that all the rightmost positive
windows in the BR correspond to either the pair $(2,0)$ and belong to
the $n^{\mbox{th}}$ column (Case I), or to the pair $(2, r)$ where
$r>0$ (Case II). Claim~\ref{clm:nineth} indicates that all positive
windows of the last column of the BR are connected (with some paths)
to each other.

We begin our considerations for Case I, first. Property 3(b) and
Claim~\ref{clm:nineth} imply that for every window $w_{i,j}^r$ from
BR, with the exception of the windows from the last (the $n^{\mbox{th}}$)
column (that is, from the $(rn+j)^{\mbox{th}}$ column, where $rn+j
<n$), there exists a path which connects it with a window from the
next  (the $(rn+j+1)^{\mbox{th}}$) column. Consequently, every positive
window in the BR is connected to a window from the last column.  This
implies that all these windows belong to the same component.

As for Case II, let $rn+k$, where $r>0$ and $1 \leq k \leq n$, be the
ordinal number of the last column of the BR (whose windows are all
assigned to the pair $(2,r)$).  Property 3(f) implies the existence of
the positive window from the BR in the previous
($(rn+k-1)^{\mbox{th}}$) column which is connected to the positive
windows from the $(rn+k)^{\mbox{th}}$ column.  If all the windows from
this previous column are assigned to the pair $(2,r)$, then we use
Claim~\ref{clm:nineth}. Otherwise we use Property 3(b) to conclude
that all the windows from this column are connected to the windows of
the last column.  Now, we can obtain the same conclusion by using the
$(rn+k-2)^{\mbox{th}}$ and $(rn+k-1)^{\mbox{th}}$ columns instead of
the $(rn+k-1)^{\mbox{th}}$ and $(rn+k)^{\mbox{th}}$ columns.  We
continue this procedure till we reach the rectangle ${\cal R}_0$, that
is, finishing with the $n^{\mbox{th}}$ and $(n+1)^{st}$ columns.  This
way all the windows from the $n^{\mbox{th}}$ column and to right of it
are connected. For the rest of the windows (from the columns to the
left of the $n^{\mbox{th}}$ column) we use Property 3(b) and
Claim~\ref{clm:nineth} similarly as in Case I.  Consequently, the
subgraph of $W_{m,n}$ induced by $B^{c,Int}$'s positive entries has a
tree structure.
\qed

\vspace*{0.5cm}

\begin{figure}[htb]
\begin{center}
\includegraphics{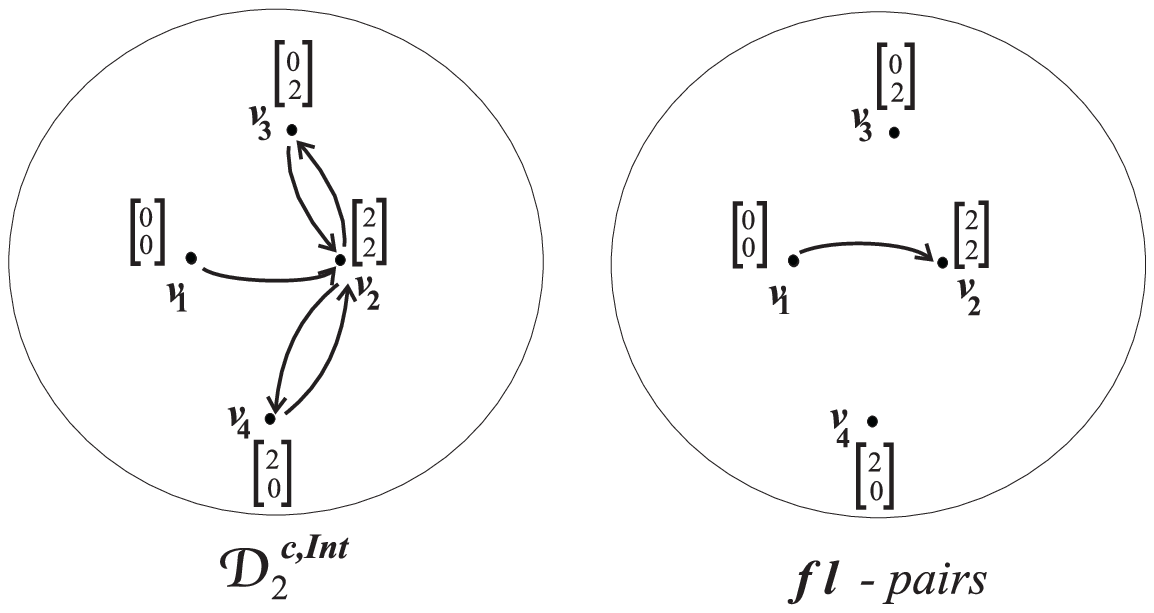}
\end{center}
\caption{The digraph $\cD^{ \; c,Int}_2$ and the corresponding set of
  pairs $\cF_2 \cL_2$.}
\label{fig:caseIm=2}
\end{figure}

Let $\cF_m $ denote the set of all possible first columns of $B^{ \;
c,Int}$.  We already know that the set of all possible columns of
the matrix $B^{ \; c,Int}$ (as defined by Properties 1--4 or 1--5
above) forms the vertex set of a digraph $\cD^{ \; c,Int}_m$.  Note
that $\cF_m \subseteq V(\cD^{ \; c,Int}_m)$.  Furthermore, the
directed edges are determined by the adjacency conditions.  Let $\cF
\cL_m$ denote the subset of $V(\cD^{ \; c,Int}_m) \times V(\cD^{ \;
c,Int}_m)$ consisting of all possible pairs of first and last
columns (which we call the \textbf{\emph{fl-pairs}}) of $B^{ \;
c,Int}$ determined by the specific properties of the first and last
columns.  The enumeration of HC$^{\; c}$'s on $P_{m+1}\times C_n$
basically comes down to the enumeration of oriented walks of length
$n-1$ in the digraph $\cD^{ \; c,Int}_m$ with the initial and last
vertices from the set $\cF \cL_m$.  For $m=2$, see
Figure~\ref{fig:caseIm=2}.  Finally, this number $\varphi_m^{ \;
c,Int}(n-1)$ should be multiplied by $n$ to obtain the correct
number of HC$^{\; c}$.

Note that the size of $\cD^{ \; c,Int}_m$ depends on whether we have
imposed the additional conditions from Property~5 on the vertices and
edges. The previously mentioned properties are quite handy,
particularly when it comes to generating the set of vertices of $\cD^{
\; c,Int}_m$.  Owing to them, it is possible to reduce the number of
edges in the said diagraph. In other words, we are actually able to
exclude the superfluous edges.


\vspace{1.0cc}
\begin{center}
{\bf 4.  THE NUMBER OF Color$^r$ WORDS OF FIXED LENGTH \\
AND CATALAN NUMBERS}
\end{center}

A color word was defined in \cite{BKP2} as a word of length $k$ over
the alphabet $\{ 2, 3, \ldots , k+1 \}$ with the following properties:
  \begin{itemize}
  \item \emph{P1}:
    If the letter $s\geq3$ appears in a word, then each
    letter from the set $\{2,3,\ldots,s-1\}$ must appear at least once
    prior to the first occurrence of $s$.  Consequently, if $i_s$
    denotes the position at which the first occurrence of $s$ can be
    found, then we must have $i_2<i_3<i_4<\cdots\,$.
  \item \emph{P2}:
    If $abab$ is a subword of a word, then $a=b$.  In every other word
    in which $a\neq b$, $abab$ cannot appear as subword.
  \end{itemize}
The number of color words of length $k$ is determined by  the
$k^{\mbox{th}}$ Catalan number $C_k = \frac{1}{k+1} \binom{2k}{k}$
\cite{BKP2}.

In \cite{BKDP} and in the previous sections we have introduced the
notions of a {\em positive (resp., negative)} or just {\em truncated
word} and a {\em positive (resp., negative)} or just {\it color$^r$
word}.  These depend on the type of HC in question (HC$^{\; nc}$ or
HC$^{\; c}$) as well as on the type of coding applied (HC$^{\; nc}$,
HC$^{\; c, Ext}, $HC$^{\; c, Int}$).  Throughout the whole process of
generating the vertices of $\cD^{ \; c,Ext}_m$ and $\cD^{ \;
  c,Int}_m$, as well as of $\cD^{nc}_m$ (which was described in PART
I) we need to construct the set of all (positive/negativ/-) color$^r$
words of length $k$. Now, we want to find the upper bound of this
set's cardinality.

In case of $B^{ \; c,Int}_m$ (for all possible $r$), this set is
determined by P1 (in accordance with Property~2 (a)) and P2 (in
accordance with Property~5(a)).  As a result, the upper bound of this
set's cardinality is precisely $C_k$.

\begin{prp}
The upper bound of the cardinality of color$^r$ words of length $k$ in
case of $B^{ \; c,Int}_m$ (for all possible $r$) is the
$k^{\mbox{th}}$ Catalan number $C_k = \frac{1}{k+1} \binom{2k}{k}$.
\end{prp}

Exactly the same situation occurs in $B^{nc}$ and $B^{ \; c,Ext}_m$
when $r<0$.  In the latter one we use the term ``positive'' or
``negative color$^r$ words'' in place of the ``color$^r$ word'' term.
This is in accordance with Properties~2 (a) and~5 (a) for both matrix
$B^{nc}$ and $B^{ \; c,Ext}_m$.

If $r=0$ and a color$^0$ word (or a positive/negative color$^0$ word)
of length $k$ is not assigned to the last ($n^{\mbox{th}}$) column of
$B^{nc}_m$ ($B^{ \; c,Ext}_m$), then this word is a word of length $k$
over the alphabet $\{ 1,2, 3, \ldots , k+1 \}$ having an additional
property, apart from P1 and P2:
  \begin{itemize}
  \item \emph{P3}:
    If $a1a$ is a subword of the word of length  $k$, then $a=1$.
  \end{itemize}
This is in accordance with Property~5 (b) of the matrix $B^{nc}$ ($B^{
  \; c,Ext}_m$).  If we add $1$ in front of each considered word, then
P3 can be interpreted as P2, but for an augmented alphabet. Therefore,
the upper bound of the cardinality of color$^0$ words of length $k$ in
these cases is $C_{k+1}$.

Last but not least, if $r>0$ (or $r=0$ and the word is assigned to the
last ($n^{\mbox{th}}$) column of $B^{nc}_m$ ($B^{ \; c,Ext}_m$)), then
the set of all color$^r$ words (positive color$^r$ words) of length
$k$ can be described as the subset of the set of all words from the
alphabet $\{1,2,\ldots,k+1\}$ that contain at least one letter $1$,
and satisfy P1, P2, and P3. This is in accordance with Properties 2
(c) and 4 (d) of the matrix $B^{nc}$ (or in accordance with Properties
2 (c) and 4 (e) of the matrix $B^{ \; c,Ext}_m$).  Note that the
number of all the negative color$^r$ words of length $k$ is equal to
the number of all the positive color$^r$ words of length $k$.
Therefore, from the previous two cases, we determine that the upper
bound of these color$^r$ words' cardinality in this case is $ C_{k+1}
- C_k$.  This way we have proved the following:

\begin{prp}
The upper bound of the cardinality of color$^r$ words (positive or
negative color$^r$ words)  of length $k$ belonging to the
$l^{\mbox{th}}$ column of  $B^{nc}$ ($B^{ \; c,Ext}_m$) is
  \[  \left\{\begin{array}{ll}
  C_k        & \mbox{ if } r<0,  \\
  C_{k+1}     & \mbox{ if } r=0 \mbox{ and } l \neq n, \\
  C_{k+1}-C_k & \mbox{ if } r>0 \mbox{ or }  (r=0 \mbox{ and } l = n).
  \end{array}\right. \]
\end{prp}

\bigskip
The words which satisfy~P1 and~P2 are called the non-interlocking and
non-skipping columns in \cite{QK}; the interpretation of
$C_{k+1}-C_k$, in the same paper, provides an alternative proof for
the case of $r>0$.


\vspace{1.0cc}
\begin{center}
{\bf 5. COMPUTATIONAL RESULTS}
\end{center}

The technique we use to compute $\ds \cH^c_m(x) \stackrel{\rm def}{=}
 \sum_{n\geq1} h^c_m(n+1) x^n$, the generating function for the
contractible HC's is, technically speaking, essentially the same as
the one utilized in Part I. For that reason, we shall only discuss a
few dissimilarities here, from the data obtained through the use of a
computer.

The primary goal of Topological Properties is to shorten the search
process throughout the digraph.  Note that they are, in fact, not
necessary for the determination of HC$^{\; c}$'s or HC$^{\; nc}$'s.
However, their importance role is to reduce the digraph's dimension
to a reasonable size by eliminating all the irrelevant vertices and
edges that cannot occur in generating any HC$^{\; c}$.

Based on all of the above theory and considerations, we wrote computer
programs to generate the matrices $M^{ \; c,Ext}_m$ and $M^{\;c,Int}_m$,
together with the adjacency matrices of the digraphs $\cD^{\;c,Ext}_m$
and $\cD^{ \; c,Int}_m$.  The dimensions of $\cD^{\;c,Ext}_m$ and
$\cD^{\;c,Int}_m$ are collected in Tables~\ref{table1}, for some
reasonable values of $m$.

The computation was performed on a personal computer equipped
with an Intel(R) Core (TM) i7-4712MQ processor (running at a speed
of 2.30GHz) with 6.00 GB of RAM, and run on a 64-bit operating
system.

Similar to  the case of the HC$^{\; nc}$'s, for the HC$^{\; c}$'s by
coding the interior tree we find that the $\cF_m = \cF^{ \; c,Int}_m
\subseteq V(\cD_m^{ \; c,Int})$.  However, when coding the exterior
trees, we came to realise that $\cF^{ \; c,Ext}_m \cap V(\cD_m^{nc}) =
\emptyset$.  The reason behind it is that the first row of the matrix
$[b_{i,j}]_{m\times n}$ has only one positive number which must be the
entry $b_{1,1}$.

\begin{table}[htb]
\begin{center}
\begin{tabular}{|c||*{8}{r|}}  \hline
$m$ & \multicolumn{1}{c|}{2} & \multicolumn{1}{c|}{3}
    & \multicolumn{1}{c|}{4} & \multicolumn{1}{c|}{5}
    & \multicolumn{1}{c|}{6} & \multicolumn{1}{c|}{7}
    & \multicolumn{1}{c|}{8} & \multicolumn{1}{c|}{9} \\ \hline\hline
    $|V(\cD^{ \; c,Ext}_m)|$
  & 3 & 11 & 44 & 174 & 644 & 2488 & - & - \\ \hline
  $|\cF^{ \; c,Ext}_m| $ & 1 & 3 & 7 & 28  & 92 & 341 & - & -  \\ \hline
$|E(\cD^{ \; c,Ext}_m)|$
  & 4  & 24 & 123 & 677 & 3446 & 18569 & - & - \\ \hline
$|\cL \cF \cS_m|$
  & 1 & 12  & 49 & 406 & 2461  & 19913  & - & - \\ \hline \hline
$|V(\cD^{ \; c,Int}_m)|$
  & 4 & 10 & 33 & 104 & 318 & 985 & 3121 & 9943 \\ \hline
$|E(\cD^{ \; c,Int}_m)|$
  & 5 & 23 & 96 & 423 & 1792 & 7857 & 34505 & 153500 \\ \hline
$|\cF \cL_m|$
  & 1 & 6 & 18 & 80 & 325 & 1413 & 6083 & 26583 \\ \hline
\end{tabular}
\end{center}
\caption{The characteristics of digraphs $\cD^{ \; c,Ext}_m$ and $\cD^{ \; c,Int}_m$.}
\label{table1}
\end{table}

Our findings for $m\leq4$ were confirmed by manual computations.  The
results displayed below agree with the values of $h^{ \; c}_m(n)$ for
$m\leq9$ and $n\leq10$ obtained in \cite{BKP}, as well as with the
values $h_m(n) = h^{nc}_m(n) + h^{ \; c}_m(n)$ for $m \leq 9$ and $11
\leq n \leq 22$ in \cite{Ka}.

Recall that the number $\varphi^c_m(k)$ represents the number of
HC$^{\; c}$'s in $P_{m+1} \times C_{k+2}$ with $w_{11}$ as the up root
of the split tree.  But, in the case of coding the exterior trees, it
represents the number of oriented walks of length $k$ in the digraph
$\cD^{c,Ext}_m$ with the pairs of initial and final vertices which are
respectively the third and first coordinates of the special
triples. Similarly, when coding the interior tree, it represents the
number of oriented walks of length $k+1$ in the digraph $\cD^{ \;
  c,Int}_m$ with the initial and last vertices from some special sets.
Hence, we label the coefficients $\varphi^c_m(k)$, where $k\geq 0$, of
the generating function $\ds \Phi^c_m(x) \stackrel{\rm def}{=}
\sum_{k\geq 0} \varphi^c_m(k)x^k$ with $\varphi^{\; c, Ext}_m(k)$ in
the coding by exterior trees, and with $\varphi^{\; c, Int}_m(k+1)$ in
the coding with the interior tree.  Generating the digraphs $\cD^{ \;
  c,Int}_7$, $\cD^{ \; c,Int}_8$, and $\cD^{ \; c,Int}_9$ requires 7
seconds, 2 minutes, and 39 minutes, respectively.


\vspace{0.5cc}
\begin{center}
{\bf 5.1. Thick Cylinder \boldmath{$P_2 \times C_n$  ($m=1$)}}
\end{center}

For $m=1$, it is easy to show that $h^c_1(n)=n$ for all $n\geq1$.
Since
  \[  h^{nc}_1(n) = \left\{\begin{array}{ll}
  2 & \mbox{if $n$ is even}, \\
  0 & \mbox{if $n$ is odd},
  \end{array}\right. \]
we can write $h_1(n) = n + 1 + (-1)^n$.


\vspace{0.5cc}
\begin{center}
{\bf 5.2. Thick Cylinder \boldmath{$P_3 \times C_n$ ($m=2$)}}
\end{center}

The digraph $\cD^{\; c,Ext}_2$ is displayed in Figure~\ref{fig:caseEm=2}.

The incidence matrix $M^{ \; c,Ext}_2$ of the corresponding digraph
$\cD^{\; c,Ext}_2$ is of order~3. The set of all possible triplets
$(l,f,s)$ has only one element $(v_1,f_1,v_1)$.  Using a similar
technique as in the case of the NC-type of HC's, we obtain
  \begin{equation}
  h^c_2(n) = \frac{n}{4} {\sqrt{2}\:}^n [1+(-1)^n].
  \label{eqn:m2c}
  \end{equation}
Since   $ h^{nc}_2(n) = 2^n - 2 $ \cite{BKDP} (Part I), from
(\ref{eqn:m2c}), we determine that
  \[ h_2(n) = 2^n - 2 + \frac{n}{4} {\sqrt{2}\:}^n [1+(-1)^n] \]
for all integers $n\geq1$.  Identical results can be obtained from
$M^{ \; c,Int}_2$ (see Figure~\ref{fig:caseIm=2}).

\vspace{0.5cc}
\begin{center}
{\bf 5.3. Thick Cylinder \boldmath{$P_4 \times C_n$ ($m=3$)}}
\end{center}

In this subsection, we provide a detailed discussion for the case of
$m=3$.  We study the HC's of type C, with coding carried out on the
exterior region first; and then we move on to the coding of the
interior region.

\bigskip\noindent
\underline{{\bf Coding the Exterior Region}}

\bigskip
We find $V(\cD^{ \; c,Ext}) = \{v_1,v_2,\ldots,v_{11}\}$; the vertices
and the adjacency matrix $M_3^{ \; c,Ext} = [m^{ \; c,Ext}_{i,j}]$ are
listed below:
\[ \begin{array}{l}
v_1 = (0^0,0^0,0^0) \\
v_2 = (0^0,1^0,0^0) \\
v_3 = (0^0,1^0,1^0) \\
v_4 = (0^0,2^0,0^0) \\
v_5 = (0^0,1^1,0^0) \\
v_6 = (-1^0,-1^0,0^0) \\
v_7 = (-1^0,0^0,1^0) \\
v_8 = (0^0,-2^{-1},0^0) \\
v_9 = (0^0,-2^0,0^0) \\
v_{10} = (0^0,2^{-1},0^0) \\
v_{11} = (0^0,-1^0,0^0)
\end{array}
\qquad \left[\begin{array}{*{11}{c}}
0&0&1&1&0&1&1&1&1&1&0 \\
1&1&0&0&0&0&0&0&0&0&0 \\
1&1&0&0&0&0&0&0&0&0&0 \\
0&0&1&1&0&0&0&0&0&0&0 \\
1&0&0&0&1&0&0&0&0&0&0 \\
1&0&0&0&0&0&0&0&0&0&1 \\
1&0&0&0&0&0&0&0&0&0&0 \\
0&0&0&0&0&0&0&1&0&0&0 \\
0&0&0&0&0&1&0&0&1&0&0 \\
0&0&0&0&1&0&0&0&0&1&0 \\
1&0&0&0&0&0&0&0&0&0&1
\end{array}\right]. \]
We also find $\cF_3 = \{f_1,f_2,f_3\}$, $|\cL\cF \cS_3| = 12$,
  \begin{eqnarray*}
  \cL \cF \cS_3
  &=&
  \{(v_{10},f_1,v_1), (v_{10},f_1,v_2), (v_1,f_1,v_1), (v_1,f_1,v_2), \\
  && \quad (v_{10},f_2,v_4), (v_{10},f_2,v_3), (v_1,f_2,v_4),
  (v_1,f_2,v_3), \\
  && \quad (v_3,f_3,v_1), (v_3,f_3,v_5), (v_2,f_3,v_1),
  (v_2,f_3,v_5)\},
  \end{eqnarray*}
where $f_1 = (1^0,1^0,1^0)$, $f_2 = (2^0,2^0,0^0)$, $f_3 =
(1^1,1^1,0^0)$. The characteristic polynomial of  $M_3^{ \; c, Ext}$
is
  \[ P_3^{ \; c,Ext}(x)
  = -x^2+7x^3-22x^4+38x^5-34x^6+6x^7+18x^8-18x^9+7x^{10}-x^{11}. \]
This implies that the sequence $\varphi^{ \; c,Ext}_3(n)$ satisfies a
recurrence relation of order~9.

\vspace{0.5cc}\noindent
\underline{{\bf Coding the Interior  Region}}

\bigskip
We find $V(\cD^{ \; c,Int}) = \{v_1,v_2,\ldots,v_{10}\}$; the vertices
and the adjacency matrix $M_3^{ \; c,Int} = [m^{ \; c,Int}_{i,j}]$ are
listed below:
\eject  
\[ \begin{array}{l}
v_1 = (0^0,0^0,0^0) \\
v_2 = (0^0,0^0,2^0) \\
v_3 = (0^0,0^0,2^1) \\
v_4 = (2^0,0^0,2^{-1}) \\
v_5 = (2^0,0^0,3^0) \\
v_6 = (2^0,2^0,2^0) \\
v_7 = (0^0,2^0,0^0) \\
v_8 = (2^0,0^0,0^0) \\
v_9 = (2^0,0^0,2^0) \\
v_{10} = (2^0,0^0,2^1)
\end{array}
\qquad
\left[\begin{array}{*{10}{c}}
 0&0&0&1&1&1&0&0&0&0 \\
 0&0&0&0&1&1&0&0&0&0 \\
 0&0&0&0&0&0&0&1&0&1 \\
 0&0&0&1&0&0&0&0&0&0 \\
 0&0&0&0&1&1&0&0&0&0 \\
 0&1&0&0&0&0&1&1&1&0 \\
 0&0&0&0&0&1&0&0&0&0 \\
 0&0&0&1&1&1&0&0&0&0 \\
 0&1&0&0&0&0&0&1&1&0 \\
 0&0&0&0&0&0&0&1&0&1
\end{array}\right] \]
We also find $\cF_3 \cL_3
  = \{(v_1,v_6),(v_1,v_9),(v_2,v_4),(v_2,v_8),(v_3,v_6),(v_3,v_9)\}. $
The characteristic polynomial for $M_3^{ \; c,Int}$ is
  \[ P_3^{ \; c,Int}(x)
  = -x^4 (1-x)^2 (1-2x+2x^2+2x^3-x^4). \]

\vspace{0.5cc}\noindent
\underline{{\bf Common Results for Both Types of Coding}}

\bigskip
Recall that $h^c_3(n) = n \varphi ^{ \; c,Ext}_3(n-2) =  n \varphi ^{
\; c,Int}_3(n-1)$,
  \[ \Phi^{ \; c,Ext}_3(x)
  = \sum_{n\geq0} \varphi_3^{ \; c,Ext}(n) x^n,
  \qquad\mbox{and}\qquad
  \Phi^{ \; c,Int}_3(x)
  = \sum_{n\geq0} \varphi_3^{ \; c,Int}(n) x^n. \]
The first few nonzero values are listed below.
\[
\begin{array}{l}
h^c_3(2) = 2\cdot \varphi^{ \; c,Ext}_3(0)
  = 2\cdot \varphi^{ \; c,Int}_3(1) = 2\cdot 1 = 2  \\
h^c_3(3) = 3\cdot \varphi^{ \; c,Ext}_3(1)
  = 3\cdot \varphi^{ \; c,Int}_3(2) =3\cdot 4 = 12 \\
h^c_3(4) = 4\cdot \varphi^{ \; c,Ext}_3(2)
  = 4\cdot \varphi^{ \; c,Int}_3(3) = 4\cdot 12 = 48 \\
h^c_3(5) = 5\cdot \varphi^{ \; c,Ext}_3(3)
  = 5\cdot \varphi^{ \; c,Int}_3(4) = 5\cdot 32 = 160 \\
h^c_3(6) = 6\cdot \varphi^{ \; c,Ext}_3(4)
  = 6\cdot \varphi^{ \; c,Int}_3(5) = 6\cdot 83 = 498 \\
h^c_3(7) = 7\cdot \varphi^{ \; c,Ext}_3(5)
  = 7\cdot \varphi^{ \; c,Int}_3(6) = 7\cdot 212 = 1484 \\
h^c_3(8) = 8\cdot \varphi^{ \; c,Ext}_3(6)
  = 8\cdot \varphi^{ \; c,Int}_3(7) = 8\cdot 540 = 4320 \\
h^c_3(9) = 9\cdot \varphi^{ \; c,Ext}_3(7)
  = 9\cdot \varphi^{ \; c,Int}_3(8) = 9\cdot 1372 = 12348 \\
h^c_3(10) = 10\cdot \varphi^{ \; c,Ext}_3(8)
  = 10\cdot \varphi^{ \; c,Int}_3(9) = 10\cdot 3485 = 34850 \\
h^c_3(11) = 11\cdot \varphi^{ \; c,Ext}_3(9)
  = 11\cdot \varphi^{ \; c,Int}_3(10) = 11\cdot 8848 = 97328 \\
h^c_3(12) = 12\cdot \varphi^{ \; c,Ext}_3(10)
  = 12\cdot \varphi^{ \; c,Int}_3(11) = 12\cdot 22464 = 269568
\end{array} \]
The generating functions are  obtained for both cases in the usual
way:
  \[ \Phi^{ \; c,Int}_3(x)
  = x \,\Phi^{ \; c,Ext}_3(x)
  = \frac{x(1+x)}{(1-x)(1-2x-2x^2+2x^3-x^4)}. \]
Then
  \begin{eqnarray*}
  \lefteqn{
  \cH_3^c(x)
  = \frac{d}{dx} \left(x^2\,\Phi^{ \; c,Ext}_3(x)\right)
  = \frac{d}{dx} \left(x\,\Phi^{ \; c,Int}_3(x)\right)} \\
  && \quad {}=
  \frac{2x(1-3x^2-2x^3+3x^4-x^6)}{(1-x)^2(1-2x-2x^2+2x^3-x^4)^2}.
  \end{eqnarray*}
From (3) in \cite{BKDP} and  the above equality  we obtain
  \begin{eqnarray*}
  \cH_3(x)
  &=& \cH_3^{nc}(x) + \cH_3^c(x) \\
  &=& 2x(2-6x+3x^2-36x^3+97x^4+96x^5-372x^6 \\
  &&{} +96x^7+280x^8-142x^9+64x^{10} -252x^{11}+132x^{12} \\
  &&{} +168x^{13}-193x^{14}+64 x^{15}-11 x^{16}+6x^{17}-2x^{18})\;/ \\
  &&{} [(1-x)^2 (1+x-x^2) (1-x-x^2) (1-2x-2x^2+2x^3-x^4)^2 \\
  &&{} (1-x-3x^2-x^3+x^4) (1+x-3x^2+x^3+x^4)].
  \end{eqnarray*}
Its power series expansion is
  \[ \cH_3(x)
  = 4x+24x^2+306x^3+850x^4+7010x^5+18452x^6+126426x^7+\cdots\,. \]


\vspace{0.5cc}
\begin{center}
{\bf 5.4. Thick Cylinder \boldmath{$P_5\times C_n$ ($m=4$)}}
\end{center}

For Hamiltonian cycles of type C, the degrees of the characteristic
polynomials for $M_4^{ \; c,Ext}$ and $M_4^{ \; c,Int}$ are~44 and~33,
respectively; they determine recursions of order~28 and~16,
respectively.  However, their generating functions $\Phi^{ \;
  c,Ext}_4(x)$ and $\Phi^{ \; c,Int}_4(x)$ indicate the same recursion
of order~12 (which was expected) for the sequences $\varphi^{ \;
c,Ext}_4(n)$ and $\varphi^{ \; c,Int}_4(n)$. Since $\varphi^{ \;
c,Ext}_4(n-2) = \varphi^{ \; c,Int}_4(n-1)=h^{ \; c}_4(n)/n$, it is
clear that
  \[ \cH_4^c(x)
  = \frac{d}{dx} \left( x^2 \,\Phi^{ \; c, Ext}_4(x)\right)
  = \frac{d}{dx} \left( x \,\Phi^{ \; c, Int}_4(x)\right), \]
where
  \[ \Phi_4^c(x) = \frac{x(1+16x^2-48x^4-8x^6+77x^8-8x^{10}+2 x^{12})}
  {(1-x) (1+x) (1-3x^2)^2 (1-11x^2-2x^6)}. \]
Thus,
  {\small
  \begin{eqnarray*}
  \cH_4^c(x)
  &=& 2x(1+35x^2-419x^4+791x^6+1251x^8-6807x^{10}+9747x^{12} \\
  &&{} -5055x^{14}+1032x^{16}+168x^{18}+36x^{20}-12x^{22}) \;/ \\
  &&{} [(1-x)^2 (1+x)^2 (1-3x^2)^3 (1-11x^2-2x^6)^2] \\
  &=& 2x+136x^3+2832x^5+44288x^7+621720x^9+8268432x^{11} \\
  &&{} +106467592x^{13}+1341213504x^{15}+16625223000x^{17} \\
  &&{} +203511990480x^{19}+2466221656712x^{21}+29639129297760x^{23} \\
  &&{} +353729229308728x^{25}+4196610165544048x^{27} \\
  &&{} +49534151824335720x^{29}+ \cdots\,.
  \end{eqnarray*}}%
From $\cH_4(x) = \cH_4^{nc}(x) + \cH_4^c(x)$, we determine
\eject  
  {\small
  \begin{eqnarray*}
  \cH_4(x)
  &=& 2x(2-4x-x^2-399x^3+264x^4+12003x^5-17018x^6-132589x^7 \\
  &&{} +242972x^8+741418x^9-1592134x^{10}-2249196x^{11}+5298156x^{12} \\
  &&{} +3895631x^{13}-6944234x^{14}-7550348x^{15}-9067499x^{16} \\
  &&{} +32222212x^{17}+44031695x^{18}-104031851x^{19}-46708066x^{20} \\
  &&{} +189953465x^{21}-27366464x^{22}-195625961x^{23}
       +111621212x^{24} \\
  &&{} +106177438x^{25}-112507558x^{26}-19954314x^{27}+54710202x^{28} \\
  &&{} -6461551x^{29}-11552268x^{30}+1712058x^{31}-156451x^{32} \\
  &&{} +2433592x^{33}-1096202x^{34}-349808x^{35}+355012x^{36}+47076x^{37} \\
  &&{} -57312x^{38}+6792x^{39}-5340x^{40}+2208x^{41}-72x^{42}+96x^{44}) \;/ \\
  &&{} [(1-x)^2 (1+x)^2 (1-2x) (1+2x) (1-2x^2) (1-3x^2)^2 \\
  &&{} (1-3x+x^3)(1-3x^2+x^3-x^4)(1-4x^2+2x^3-2x^4-x^5) \\
  &&{} (1-5x+2x^2+8x^3-8x^4+x^5-x^6) (1+11x^2+2x^6)^2].
  \end{eqnarray*}}%
Upon expansion, we obtain
  {\small
  \begin{eqnarray*}
  \lefteqn{\cH_4(x)
  = 4x+24x^2+306x^3+850x^4+7010x^5+18452x^6+126426x^7+351258x^8} \\
  &&{} +2127332x^9+6355404x^{10}+35085590x^{11}+112481980x^{12}
       +577875650x^{13} \\
  &&{} +1970896234x^{14}+9576146794x^{15}+34373120896x^{16}
       +160047128522x^{17} \\
  &&{} +598167523522x^{18}+2697774177200x^{19}+10398653965136x^{20} \\
  &&{} +45813998934398x^{21}+180683364527008x^{22}
       +782729112571558x^{23} \\
  &&{} +3138757868554550x^{24}+13435232382112114x^{25}
       +54519162573345144x^{26} \\
  &&{} +231410726096158954x^{27}+946929235189639806x^{28} \\
  &&{} +3995898137059583288x^{29}+16446553366281600876x^{30}+\cdots\,.
  \end{eqnarray*}}%


\vspace{0.5cc}
\begin{center}
{\bf 5.5. Thick Cylinder \boldmath{$P_6\times C_n$ ($m=5$)}}
\end{center}

For the HC's of type C, the characteristic polynomial for $M_5^{ \;
  c,Ext}$ (of order~174) yields a recurrence of order~140, and the
characteristic polynomial for $M_5^{ \; c,Int}$ (of order~104)
determines a recurrence of order~68.  The generating functions $\Phi^{
\; c,Ext}_5(x)$ and $\Phi^{ \; c,Int}_5(x)$ indicate a recursion of
order~48.
{\small
\begin{eqnarray*}
\lefteqn{%
\Phi_5^c(x) =
x(1+6x-18x^2-414x^3+848x^4+6554x^5-16045x^6-37690x^7} \\
&&{} +103281x^8+128504x^9-265355x^{10}-672050 x^{11}+502008x^{12}
     +3340076x^{13} \\
&&{} -2448123x^{14}-9954494x^{15}+10693383x^{16}+18205338x^{17}
     -29135762x^{18} \\
&&{} -20019118 x^{19}+53993028x^{20}+8041536x^{21}-72266964x^{22}
     +14531980x^{23} \\
&&{} +71713080x^{24}-34388270x^{25}-52305506x^{26}+39831436x^{27}
     +26241113x^{28} \\
&&{} -30511612x^{29}-7246130x^{30}+16104692x^{31}-373346x^{32}
     -5790448x^{33} \\
&&{} +1221241x^{34}+1368750x^{35}-510697x^{36}-198552x^{37}+111455x^{38}
     +14804x^{39} \\
&&{} -13910x^{40}-84x^{41}+945x^{42}-72x^{43}-27x^{44}+4x^{45})\;/ \\
\noalign{\eject}
&&{} [(1-x) (1-x-x^2) (1-2x-x^2+x^3) (1+2x-x^2-x^3) \\
&&{} (1-2x-x^2+2x^3-x^4) (1-x-7x^2+2x^3-2x^4) (1-2x^2+3x^3+x^4-x^5) \\
&&{} (1-4x+x^2+6x^3-4x^4-2x^5+x^6) (1+4x+x^2-6x^3-4x^4+2x^5+x^6) \\
&&{} (1-5x-14x^2+63x^3-12x^4-90x^5+35x^6+66x^7-118x^8+8x^9 \\
&&{} +82x^{10}-42x^{11}-28x^{12}+4x^{13}-2x^{14})]
\end{eqnarray*}}%
Below is the power series expansion of $\cH_5^c(x)$.
{\small
\begin{eqnarray*}
\lefteqn{\cH_5^c(x) =
2x+48x^2+612x^3+4520x^4+35964x^5+229698x^6+1575288x^7+9806292x^8} \\
&&{} +62999960x^9 +387822094x^{10}+2411860680x^{11}+14706401372x^{12} \\
&&{} +89805227764x^{13}+542922027450x^{14}+3277207263040x^{15}
     +19667429401654x^{16} \\
&&{} +117755148280932x^{17}+702348082721928x^{18}+4179389353497440x^{19} \\
&&{} +24800669448996294x^{20}+146859912517712812x^{21}
     +867755187436181848x^{22} \\
&&{} +5117982982251905808x^{23}+30131949609066739700x^{24} \\
&&{} +177121683074273170272x^{25}+1039599437405096152836x^{26} \\
&&{} +6093526747267596471744x^{27}+35670915471779662426386x^{28} \\
&&{} +208566648331009396119000x^{29}+\cdots\,.
\end{eqnarray*}}%


\vspace{0.5cc}
\begin{center}
{\bf 5.6. Thick Cylinder \boldmath{$P_{7}\times C_n$ ($m=6$)}}
\end{center}

We obtain
{\small
\begin{eqnarray*}
\lefteqn{\cH_6^c(x) =
2x+2032x^3+263736x^5+22337664x^7+1641664580x^9+113092326312x^{11}} \\
&&{} +7512031798348x^{13}+487293888097600x^{15}
     +31078838281479156x^{17} \\
&&{} +1956749096194717760x^{19}+121942699478516467980x^{21} \\
&&{} +7535939697350674950480x^{23}+462464193503836875708188x^{25} \\
&&{} +28212097969001607154778424x^{27}
     +1712255987823212304590396640x^{29}+\cdots\,.
\end{eqnarray*}}%
which was derived from
{\small
\begin{eqnarray*}
\lefteqn{\Phi_6^c(x) =
x(1+306x^2-40690x^4+2088888x^6-59854356x^8+1041724854x^{10}} \\
&&{} -9963350575x^{12}-1561595514x^{14}+1764222372901x^{16}
     -33359541871130x^{18} \\
&&{} +391283632798625x^{20}-3409940836072834x^{22}
     +23527202411977523x^{24} \\
&&{} -132804217617691704x^{26}+626327321659400394x^{28}
     -2507190853725762016x^{30} \\
&&{} +8634200715329254103x^{32}-25906190273641652336x^{34} \\
&&{} +68559708278394067292x^{36}-161854074407081021262x^{38} \\
&&{} +343843110397151257836x^{40}-660259656542874312136x^{42} \\
&&{} +1145535722603029360938x^{44}-1788319405800757683806x^{46} \\
&&{} +2496733228293042684759x^{48}-3090833148880542271276x^{50} \\
&&{} +3332184138797783431832x^{52}-2977551413530470915268x^{54} \\
&&{} +1874479143252245895283x^{56}-135583843251370310752x^{58} \\
&&{} -1725296888982641415559x^{60}+2978301300202755230596x^{62} \\
&&{} -3196573579330611658014x^{64}+2594775826885761943406x^{66} \\
\noalign{\eject}  
&&{} -1764806452649464365958x^{68}+1132024373770684603442x^{70} \\
&&{} -747388500864816970949x^{72}+488605725204679503240x^{74} \\
&&{} -281985452668526323045x^{76}+129099929765096588046x^{78} \\
&&{} -40724401950514905533x^{80}+4766321154300357988x^{82} \\
&&{} +3107539531302455678x^{84}-2141576362869785724x^{86}
     +699992921850375378x^{88} \\
&&{} -164891856789160528x^{90}+42040968674824716x^{92}
     -11653173730620232x^{94} \\
&&{} +2342301220645576x^{96}-366020508217376x^{98}
     +106155952612144x^{100} \\
&&{} -34874926374880x^{102}+6594270112768x^{104}
     -673590339328x^{106} \\
&&{} +37441754880x^{108}+169198080x^{110}-292896768x^{112}
     +15998976x^{114})) \;/ \\
&&{} [(1-2x) (1+2 x) (1-2x^2) (1-4x^2+2x^4) (1-8x^2+14x^4) \\
&&{} (1-27x^2+225x^4-641x^6+659x^8-227x^{10}+46x^{12}-123x^{14} \\
&&{} +169x^{16}-49x^{18}+4x^{20}) (1-37x^2+397x^4-1681x^6+2639x^8 \\
&&{} -903x^{10}-56x^{12}+307x^{14}+525x^{16}-209x^{18}-8x^{20}) \\
&&{} (1-35x^2+322x^4-1485x^6+4262x^8-7682x^{10}+7755x^{12} \\
&&{} -10671x^{14}+18616x^{16}-9492x^{18}-1484x^{20}+1589x^{22}
     -62x^{24}-40x^{26}) \\
&&{} (1-85x^2+1932x^4-20403x^6+116734x^8-386724x^{10}+815141x^{12}
     -1251439x^{14} \\
&&{} +1690670x^{16}-2681994x^{18}+4008954x^{20}-3390877x^{22}
     +1036420x^{24} \\
&&{} +178842x^{26}-92790x^{28}-17732x^{30}+5972x^{32}-1728x^{34}
     -144x^{36})],
\end{eqnarray*}}%


\vspace{0.5cc}
\begin{center}
{\bf 5.7. Thick Cylinder \boldmath{$P_{m+1}\times C_n$ ($7\leq m\leq 9$)}}
\end{center}

For the sake of brevity, we only display the power series expansion of
the generating functions $\cH_m^c(x)$ for $7\leq m\leq9$.

{\small
\begin{eqnarray*}
\lefteqn{\cH_7^c(x) =
2x+192x^2+8192x^3+127860x^4+2779014x^5+35663964x^6+605992784x^7} \\
&&{} +7769376972x^8+116791523380x^9+1519170232976x^{10}
     +21412201037580x^{11} \\
&&{} +280509236582900x^{12}+3817205794180856x^{13}
     +50048772776920380x^{14} \\
&&{} +667452277157951872x^{15}+8730496956098122924x^{16} \\
&&{} +114990875591325208344x^{17}+1498721829346080971718x^{18} \\
&&{} +19577329280144309140500x^{19}+254184297263298653321994x^{20} \\
&&{} +3300736306177174727889026x^{21}
     +42700068205017640140982112x^{22} \\
&&{} +551992937500828921720192872x^{23}
     +7117443280523917273056970850x^{24} \\
&&{} +91678074802674943650656279184x^{25}
     +1178664397321648769186649515370x^{26} \\
&&{} +15136943041102404084253082680484x^{27} \\
&&{} +194109908825815965787089284625154x^{28} \\
&&{} +2486557768079418847989177095267850x^{29}+\cdots\,,
\end{eqnarray*}}%
{\small
\begin{eqnarray*}
\lefteqn{\cH_8^c(x) =
2x+29104x^3+22869384x^5+10215798448x^7+3817933082020x^9} \\
&&{} +1320093157541136x^{11}+437662770447567560x^{13}
     +141338955368771390016x^{15} \\
&&{} +44820915345596090414880x^{17}+14022558891295056887443160x^{19} \\
&&{} +4341003538573245627716733276x^{21}
     +1332438402563600063493162541728x^{23} \\
&&{} +406097228507066527913083107762828x^{25} \\
&&{} +123030579109675576558337448581510896x^{27} \\
&&{} +37082080183727206133637758977416982500x^{29}+\cdots\,,
\end{eqnarray*}}%
{\small
\begin{eqnarray*}
\lefteqn{\cH_9^c(x) =
2x+768x^2+112164x^3+3616880x^4+222067212x^5+5539931796x^6} \\
&&{} +242178636928x^7 +6169925169414x^8+224360971248960x^9 \\
&&{} +5973677282007402x^{10}+195609021230822100x^{11}
     +5368583261802264972x^{12} \\
&&{} +165442535132471644292x^{13}+4616256789338403997830x^{14} \\
&&{} +137248345001943810512192x^{15}+3858464231851630072287480x^{16} \\
&&{} +112245666066341094211887474x^{17}
     +3163164190641471563556546716x^{18} \\
&&{} +90762230874863948167367645720x^{19}
     +2557054414248684611758303990008x^{20} \\
&&{} +72707096815758305550335466105864x^{21} \\
&&{} +2045202670958705026338344754333024x^{22} \\
&&{} +57786674199846212563657479095447016x^{23} \\
&&{} +1622114630916418603623686588180361000x^{24} \\
&&{} +45620599345582502819377370379402448790 x^{25} \\
&&{} +1277757171356042779682960928336763812048x^{26} \\
&&{} +35808036061847453918421598630326003072512x^{27} \\
&&{} +1000745362534879987116390583909098103232022x^{28} \\
&&{} +27965040033048966560373497047404511628553450x^{29}+\cdots\,.
\end{eqnarray*}}%


\vspace{1.0cc}
\begin{center}
{\bf 6. ASYMPTOTIC VALUES --- A SUMMARY OF RESULTS}
\end{center}

For type~C Hamiltonian cycles, our computational data confirm that for
$2 \leq m\leq 6$, the characteristic polynomials of $M^{ \; c,Ext}$
and $M^{ \; c,Int}$ have only one (and the same) simple real positive
dominant characteristic root $\theta_{m,c}$, see Table~\ref{Table5b}
(for even $m$, there are two dominant characteristic roots
$\theta_{m,c}$ and $-\theta_{m,c}$; whereas for odd $m$, there is a
unique dominant simple characteristic root $\theta_{m,c}$).

\begin{table}[htb]
\begin{center}
\begin{tabular}{|c||c|c|c|c|} \hline
$m $ & $d.d.( \Phi^{ \; c,Int}_m) $
& $d.d.( \cH^{ \; c,Int}_m) $ & &
 \\
 & $=d.d.( \Phi^{ \; c,Ext}_m) $
& $ =d.d.( \cH^{ \; c,Ext}_m) $ &$d.d.( \cH^{ \; nc}_m) $ &$d.d.( \cH_m) $
 \\ \hline\hline
2 &  2  & 4 & 2  & 5   \\ \hline
3 & 5  &  10& 12  & 22   \\ \hline
4 & 12 & 22 & 26  & 44 \\ \hline
5 & 48 & 96  & 84  & 180 \\ \hline
6 & 114  & 228 & - & -
 \\ \hline
\end{tabular}
\end{center}
\caption{The degree of denominators ($d.d.$) of  $\Phi^{ \; c,Int}_m$
($\Phi^{ \; c,Ext}_m)$, $\cH^{ \; c,Int}_m$ ($\cH^{ \; c,Ext}_m$),
$\cH^{nc}_m$ and  $\cH_m$ for $2\leq m\leq6$.}
\label{Table5a}
\end{table}

Note that the sum of the degrees of the denominators of $\cH^{ \;
nc}_m$ and $\cH^{ \; c}_m$ does not exceed the degree of the
denominator of $\cH_m$ which determine the order of the recurrence
relation of the sequence $h_m(n)$.  This goes in favour of our
decision to split our work in two --- the problem of determining the
HC$^{\; nc}$'s and the HC$^{\; c}$'s.  However, this was not the case
for thin cylinders \cite{BKP}.

\begin{table}[htb]  \small
\begin{center}
\begin{tabular}{|c||l|l|l|} \hline
$m$&  $\theta_{m,nc}$ & $\theta_{m,c}$ \\ \hline\hline
2 & 2
  & $\sqrt{2}$ \\ \hline
3 & $2.36920540709246654628$
  & $2.53861576354917625747 $ \\ \hline
4 & $4.16748148276892815337 $
  & $3.31910824039947675342 $ \\ \hline
5 & $5.34684254175541433292 $
  & $5.65205864851675849429 $ \\ \hline
6 & $\approx_{(100)} 8.908937311$
  & $7.52634546292690578713 $ \\ \hline
7 & $\approx_{(100)} 11.8249316$
  & $\approx_{(100)} 12.382351641593$ \\ \hline
8 & $\approx_{(70)} 19.17 $
  & $\approx_{(100)}  16.77216819355$ \\ \hline
9 & $\approx_{(30)} 26$
  & $\approx_{(30)} 27$ \\ \hline
\end{tabular}
\\ \ \\ \ \\
\begin{tabular}{|c||l|} \hline
$m$ & $a_{m,c}$ \\ \hline\hline
2 & 0.25 \\ \hline
3 & $0.31357228606585772287 $ \\ \hline
4 & $0.19324623166497686532 $ \\ \hline
5 & $0.18876590435542745301 $ \\ \hline
6 & $0.14384483205795162266 $ \\ \hline
7 & $\approx_{(100)} 0.13626186172698$\\ \hline
8 & $\approx_{(100)} 0.11306933143427 $\\ \hline
9 & $\approx_{(30)}0.1053 $\\ \hline
\end{tabular}
\end{center}
\caption{The approximate values of $\theta_{m,nc}$ (with
  $a_{m,nc}=1$), $\theta_{m,c}$ and $a_{m,c}$ for $2\leq m\leq9$,
  where $\approx_{(n)}$ means the estimate based on the first n
  entries of the sequence.}
\label{Table5b}
\end{table}

The denominators of the generating functions
$\Phi_m^{ \; c,Ext}(x)$ (or $\Phi_m^{ \; c,Int}(x)$) and $\cH_m^c(x)$ have
the radius of convergence $1/\theta_{m,c}$. For $\cH_m^c(x)$, the
dominant root of the denominator has multiplicity~2 because
  \[ \cH^c_m(x)
  = \frac{d}{dx} \left(x^2 \,\Phi^{ \; c,Ext}_m(x)\right)
  = \frac{d}{dx} \left(x \,\Phi^{ \; c,Int}_m(x)\right), \]
which is deduced from
  \begin{equation}
  h^{ \; c}_m (n)
  = n \varphi^{ \; c,Int}_m (n-1)
  = n \varphi^{ \; c,Ext}_m (n-2).
  \label{eqn:a1}
  \end{equation}

Let $a_{m,c}$ denote the coefficient of $n \theta_{m,c}^n$ in the
explicit expression for $h_m^c(n)$ derived from the recurrence
relation.  From (\ref{eqn:a1}) and (\ref{thm:zeroh}), we conclude that
the coefficient for $\theta_{m,c}^n$ is 0, when $m$ is odd; and for
even $m$, the coefficients of $n(-\theta_{m,c})^n$, $\theta_{m,c}^n$,
and $(-\theta_{m,c})^n$ are $a_{m,c}$, 0, and 0, respectively.  All
this is neatly summarized below:
  \[ h^{ \; c}_m(n)
  \sim \left\{\begin{array}{ll}
  a_{m,c} n \theta_{m,c}^n & \mbox{if $m$ is odd}, \\
  2 a_{m,c} n \theta_{m,c}^n & \mbox{if $m$ is even and $n$ is even}, \\
  0 & \mbox{if $m$ is even and $n$ is odd}.
  \end{array}\right. \]

\small
\vspace{0.1cc}

\noindent

\vspace{0.5cc}
\noindent\emph{Example 8.} \ %
The number $h^{ \; c}_5(250)$ has 190 digits:
  $$  \begin{aligned}
  h^{ \; c}_5(250) =
  & \mbox{\textbf{5315308482081368176130135765364458028442269812845}}\,829751132 \\
  & 282449366037705916081244966378616480765252334858462630450547142728 \\
  & 707832260337088675894551742436677743236273632760226951744130398000
  \end{aligned} $$
as well as $a_{5,c} \cdot 250 \cdot\theta^{250}_{5,c}$, and their
first $49$ digits are identical.

\vspace{0.5cc}
\noindent\emph{Example 9.} \ %
The number $h^{ \; c}_8(100)$ has 124 digits:
  $$\begin{aligned}
  h^{ \; c}_8(100) =
  & \mbox{\textbf{650572515095530}}\,4765274909197134354116977319209669418653015912 \\
  & 096606195064869245906663377660631373746911131674517266864224600
  \end{aligned} $$
as well as $2a_{8,c} \cdot 100 \cdot\theta^{100}_{8,c}$, and their
first $15$ digits are identical.

\vspace{0.5cc}

\normalsize
With the results and conjecture of Part I \cite{BKDP} in mind,
together with the assumptions about the positive dominant
characteristic root $\theta_{m,c}$, we may now make a conjecture
regarding the behaviour of the number of all Hamiltonian cycles
$h_m(n)$ in the graph $P_{m+1} \times C_n$, when $m$ is fixed and $n
\rightarrow \infty$, as below:
  \[ h_m(n) = h^{ \; c}_m(n) + h^{ \; nc}_m (n)
  \sim \left\{\begin{array}{ll}
  a_{m,c} n \theta_{m,c}^n  + (1 + (-1)^n) \theta_{m,nc}^n
    & \mbox{if $m$ is odd} , \\  [3pt]
  a_{m,c} n (1 + (-1)^n) \theta_{m,c}^n  +  \theta_{m,nc}^n
    & \mbox{if $m$ is even }.
  \end{array}\right. \]

When $m \leq 9$ the data shows that  $\theta_{m,c}>\theta_{m,nc}$ for
odd $m$; whereas $\theta_{m,nc}>\theta_{m,c}$ for even $m$.  Assuming
that the same holds for all the values of $m$, we propose

\begin{conj}
  \[ h_m(n)= h^{ \; c}_m(n) + h^{ \; nc}_m (n)
  \sim \left\{\begin{array}{ll}
  a_{m,c} n \theta_{m,c}^n & \mbox{if $m$ is odd}, \\
  \theta_{m,nc}^n & \mbox{if $m$ is even }.
  \end{array}\right. \]
\end{conj}

\vspace{0.5cc}
\noindent\emph{Example 10.} \ %
For $m=5$ and   $n=250$, $h_5^{nc} (250)  \sim 2.1153 \cdot 10^{182}$,
whilst $h_5^{c} (250) $ $ \sim 5.3153 \cdot 10^{189}$, and so \small
$h_5(250)  = h^{ \; c}_5(250) + h^{ \; nc}_5 (250)  \sim 5.3153 \cdot
10^{189}\sim h_5^{c} (250) $.

\normalsize

\vspace{0.5cc}
\noindent\emph{Example 11.} \ %
For $m=6$ and  $n=100$, \ \small  $h_6^{nc} (100)      =
96 070 554 870 981 782 995 827 137 14212 $ \\
$8 630 238 785 786 563 044 765 223 962 050 649 940 105
158 800 796 411 036 738 881 670
\sim 9,6071 \cdot 10^{94}$,  \normalsize while \small $h_6^{c} (100)
= 131020 4649 1976349492 4652519229 6382014958 6941469996 1724074530458 $ \\
$  96850493817395 8572002437400
\sim 1,3102 \cdot 10^{89}$.
\normalsize
Thus,  $h_6(100) $ $\sim  9,6071 \cdot 10^{94} \sim  \theta_{6,nc}^{100}$.


\bigskip
From the results obtained for $m \leq 6$ we have spotted that the
positive dominant characteristic root $\theta_{m,c}$ of $h_{m}^{c}
(n)$ corresponding to $P_{m+1} \times C_n$ is the same as the positive
dominant characteristic root of the same sequence associated to
$P_{m+1} \times P_{n}$ \cite{BPPB}.  Observe that the polynomial
$1-2x-2x^2+2x^3-x^4$, being in the denominator of $\Phi^{ \;
c,Ext}_3(x)$ and $ \Phi^{ \; c,Int}_3(x)$ (or ${\cal H}^{ \;
c}_3(x)$ and ${\cal H}_3(x)$) is also in the
denominator of the generating function of sequence corresponding to
$P_4 \times P_n$ \cite{BT}. The same phenomenon occurs for $4 \leq m
\leq 6$, as well. The obtained approximate values of the dominant
characteristic root for $7 \leq m \leq 9$ speak in favour of the same
conclusion.  That brings us to our next conjecture:

\begin{conj}
Let $r_m(n)$ ($m \geq 1$) be the number of HC's in $P_{m+1} \times
P_n$ and $h_{m}^{c} (n)$ be the number of contractible  HC's in
$P_{m+1} \times C_n$ ($m \geq 1$).  Then
  \[ \lim_{n \rightarrow \infty} \frac{r_m(n)}{r_m(n-1)}
  = \lim_{n \rightarrow \infty} \frac{h_m^{c}(n)}{h_m^{c}(n-1)}
  \mbox{ for odd $m$   } \]
and
  \[ \lim_{n \rightarrow \infty} \frac{r_m(2n)}{r_m(2n-2) }
  = \lim_{n \rightarrow \infty} \frac{h_m^{c}(2n)}{h_m^{c}(2n-2)}
  \mbox{ for even $m$.  } \]
\end{conj}

If the above conjecture holds, then using merely the data acquired
from the sequence of $r_m(n)$'s we could conclude that $\theta_{10,c}
\sim 37.03764916$, $\theta_{11,c} \sim 58.75$, $\theta_{12,c} \sim
81.366569$ and $\theta_{13,c} \sim 127.7$.  In other words, we would
not require the exact value of $h_m^c(n)$ to do so.


\vspace{1.0cc}
\begin{center}
{\bf 7. CLOSING REMARKS AND FURTHER RESEARCH}
\end{center}

For the purpose of enumerating Hamiltonian cycles on $P_{m+1} \times
C_n$ we have provided one characterization of the non-contractible
HC's in Part~I, and two characterizations of the contractible HC's
with fixed up root of the split tree in $w_{11}$.

\bigskip\noindent
\underline{1. Confirmation of the old data and the process of
obtaining the new ones}

\bigskip
Both of the computer programs dealing with the HC$^c$ case have
provided the same number for $h_{m}^c(n)$, when $m \leq 7$, which
agrees with the corresponding values of \cite{BKP}. The latter holds
for $h_{m}^c(n)$ as well, where $8 \leq m \leq 9$, when obtained in
the act of coding the interior. The sum of sequences obtained in all
the three programs, i.e. the numerical values of $h_{m}(n) =
h_{m}^{nc}(n) + h_{m}^c(n)$ agree with the ones obtained earlier in
\cite{BKP}, \cite{Kar} and \cite{Ka}, for $ m \leq 9$ and $ n \leq
22$.  We have derived new data for $m \leq 9$ and $n \geq 23$.

\bigskip\noindent
\underline{2. The advantage of coding the interior over coding the
exterior}

\bigskip
Comparing the number of vertices of the digraphs $\overline{\cD}^{ \;
c,Ext}_m$ and $\cD^{ \; c,Int}_m$ (see Table 1) one can come to a
conclusion that coding the windows of the interior region is much more
efficient than the process of coding the windows of the exterior
region of a HC$^{\; c}$.

\bigskip\noindent
\underline{3. The advantage of coding the regions over coding the
vertices}

\bigskip
For the purpose of obtaining the total number of HC's, in case of thin
cylinder $C_m \times P_n$, coding the vertices has proven itself to be
a better approach. Namely, for a fixed $m$, the number of vertices of
the assigned digraph in the aforementioned approach \cite{BPPB}
turned out to be smaller than the number obtained when coding the
regions \cite{BKP}.  Additionally, the order of recursion of the total
number of HC's for thin cylinders is smaller than for special HC's,
i.e.  HC$^{\; nc}$ and HC$^{\; c}$.  The results show that the
opposite is true for thick cylinders.  This supports the choice of our
approach when tackling the thick cylinders, although it has to be
split into parts. It goes without saying, that further research in
this direction, would be nice. Particularly, it would be a good idea
to utilize the approach with coding the vertices so as to be able to
reach a precise conclusion regarding the pros and cons of coding the
regions, by a direct comparison of the number of vertices of the
assigned digraphs.

\bigskip\noindent
\underline{4. Open questions}

\bigskip
For the initial values of $m$ we have come to notice that the numbers
of HC$^{\; nc}$'s are the dominant ones for even $m$, whereas the
numbers of HC$^{\; c}$'s are such for odd $m$.  That prompted us to
make a conjecture about the asymptotic behaviour of the total number
of HC's in the graph $P_{m+1} \times C_{n}$.  Moreover, certain
matchings between the dominant characteristic roots of the sequences
the numbers of HC$^{\; c}$'s in $P_{m+1} \times C_{n}$ and $P_{m+1}
\times P_{n}$ for small values of $m$ are noticed. This way, we have
come to yet another conjecture regarding the asymptotic behaviour of
the entries of these two sequences.


\vspace{1.0cc}
\begin{center}
{\bf ACKNOWLEDGEMENT}
\end{center}
The authors would like to express their  gratitude to the referees on some useful suggestions and
helpful comments which improved the clarity of the
presentation.
The authors acknowledge financial support of the Ministry of Education,
Science and Technological Development of the Republic of Serbia
(Grant No. 451-03-9/2021-14/200125  and  451-03-68/2020-14/200156).


\vspace{2.5cc}

\noindent
\parbox[]{7cm}{
  \noindent Dept.\ of Math.\ \&\ Info., \\
  Faculty of Science,\\
  University of Novi Sad,\\
  Novi Sad, Serbia \\
  E-mail: olga.bodroza-pantic@dmi.uns.ac.rs
  \par\vspace{1cc}
  \noindent Dept.\ of Math.\ Sci.,\\
  SUNY Fredonia, \\
  Fredonia, NY 14063, U.S.A. \\
  E-mail: kwong@fredonia.edu}
 \parbox[]{5.5cm}{
  \noindent Faculty of Technical Sciences,\\
  University of Novi Sad,\\
  Novi Sad, Serbia\\
  E-mail: jelenadjokic@uns.ac.rs \\
  E-mail: rade.doroslovacki@uns.ac.rs
  \par\vspace{1cc}
  \noindent Department of Physics,\\
  Faculty of Science,\\
  University of Novi Sad,\\
  Novi Sad, Serbia \\
  E-mail: mpantic@df.uns.ac.rs}

\end{document}